\documentclass[11pt]{article}

\usepackage{bbding}
\usepackage{pifont}
\usepackage{amsfonts}
\usepackage{graphicx}
\usepackage{amsmath}
\usepackage{amsthm}
\usepackage{amssymb}
\usepackage{concrete}

\def\beq{\begin{equation}}
\def\eeq{\end{equation}}

\def\nd{\noindent}

\def\<{\leq}
\def\>{\geq}

\newtheorem{thm}{Theorem}[section]
\newtheorem{lem}{Lemma}[section]
\newtheorem{prop}{Proposition}[section]

\newtheorem{defi}{Definition}[section]
\newtheorem{rem}{Remark}[section]

\begin{document}
\title{BIRMAN-MURAKAMI-WENZL TYPE ALGEBRAS FOR GENERAL COXETER GROUPS }
\author{Zhi Chen}
\date{}
\maketitle

\begin{abstract}
\noindent We introduce a BMW type algebra for every Coxeter group,
 These new algebras are introduced as
deformations of the Brauer type algebras introduced by the author,
they  have the corresponding Hecke algebras as quotients.
\end{abstract}
\section{introduction}
The Birman-Murakami-Wenzl algebras $B_n (\tau ,l)$ (BMW
algebras)introduced by Birman,Wenzl in \cite{BW}  are natural
deformations of Brauer's centralizer algebras $B_n (\tau)$. These
two kinds of algebras were studied extensively during last two
decades. For an explanation of their backgrounds and their relation
to three dimensional Topology, consult\cite{BW} \cite{Mu}.

In 2001, H\"{a}ring-Oldenberg \cite{Ha} introduced the Cyclotomic
BMW algebras and Cyclotomic Brauer algebras associated with
$G(m,1,n)$ type pseudo reflection groups.  The Cyclotomic BMW
algebras were studied intensively in recent years. In \cite{RX} Rui
and Xu classified irreducible representations of these algebras. In
\cite{GM} Goodman proved they support cellular structures. In
\cite{Yu} Yu studied freeness of them for certain parameter ring,
and also proved cellularity of these algebras. In \cite{GH} Goodman
and Hauschild introduced Affine BMW algebras from topological
considerations and give them a presentation. In 2005, Cohen,
Gijsbers and Wales \cite{CGW1} introduced a BMW type algebra and a
Brauer type algebra for every simply laced Coxeter group. Then they
calculated dimension, and established the semisimplicity and
cellularity of these algebras (of finite type) in \cite{CGW2},
\cite{CW2}.  In this paper we present a BMW type algebra for every
Coxeter group, as the deformations of the generalized Brauer type
algebras constructed in our last paper \cite{CH}. When $W$ is of
simply laced type, our algebras are isomorphic to the algebras from
\cite{CGW1}(Proposition 5.2).

 Since the present paper is very closely related to the
simply laced BMW algebras, we write down the presentation of them
according to \cite{CGW1}. Let $\Gamma$ be a simply laced Dynkin
diagram.  Here 'simply laced' means containing no multiple bonds.
For two nodes $i,j$ of $\Gamma$, if $i$ and $j$ are connected by a
bond, we write $i\sim j$, otherwise we write $i\nsim j$.
\begin{defi}[\cite{CGW1}]

Let $\Gamma$ be a rank $n$ simply laced Dynkin diagram. The type
$\Gamma$ Brauer algebra $B_{\Gamma} (\tau )$ over $\mathbb{Q}
(\tau)$ and  the type $\Gamma$ BMW algebra $B_{\Gamma}$ over
$\mathbb{Q}(l,x)$ are defined in the following table. $$TABLE\ 1.\
Presentation\ for\ B_{\Gamma } (\tau )\  and\  B _{\Gamma }
(\tau,l).$$
\end{defi}

\begin{tabular}{|l|l|l|}

\hline

   & $B_{\Gamma } (\tau )$  & $B _{\Gamma } (\tau,l)$ \\
\hline

  Generators    & $s_i  (i\in I ) $; $e_i (i\in I) $& $X_i (i\in I)$ ; $E_i (i\in I)$ \\
\hline

  Relations   &  $s_i s_j s_i =s_j s_i s_j$, $if\ i\sim j $  ;  & $X_i X_j X_{i} =X_j X_{i} X_j$, $if\ i\sim j$;
\\
 & $s_i s_j e_i = e_j s_i s_j $ if $i \sim j$; & $X_i X_j E_i = E_j X_i X_j
 $ if $i \sim j$;\\
 & $s_i s_j =s_j s_i$ ,$if\ i\nsim j $;& $X_i X_j =X_j X_i ,if\ i\nsim j $ ;\\
 &$s_i ^2 =1 ,for\  all\  i $;& $l(X_i ^2 +\nu X_i -1) =\nu E_i , for\  all\
i $;\\
 & $s_i e_i =e_i ,for\ all\ i $;&$ X_i E_i = l^{-1} E_i ,for\ all\ i $;\\
& $e_i s_{j} e_{i} = e_{i}, if\ i\sim j $; &$ E_i X_{j} E_{i} = l
E_i ,if i\sim j $;
\\

&$s_i e_j = e_j s_i , if\ i\nsim j $; & $X_i E_j =E_j X_i ,if i\nsim j $;\\
&$e_i ^2 = \tau e_i ,for\  all\  i.$ &$ E_i ^2 =\tau E_i $
. \\

  \hline
\end{tabular}\\

  Where $\nu =\frac{l-l^{-1}}{1-\tau }$. When the Coxeter
group is of $A_{n-1}$ type,that is, being the n-th symmetric group,
the corresponding simply laced BMW algebra and Brauer algebra are
just $B_n (\tau ,l)$ and $B_n (\tau)$ respectively. In \cite{CH} the
author introduced a Brauer type algebra(see Definition 2.1,
Definition 2.2 of the present paper) for every Coxeter group and
every pseudo reflection group as candidates for generalized Brauer
type algebras. In that construction, certain KZ connections and
generalized Lawrence-Krammer representations for finite type Coxeter
groups introduced by Marin ( including a slightly further
generalization by the author) play important roles. When $W$ is a
finite pseudo reflection group, the algebra $Br_W (\Pi )$ supports a
flat, $W-$invariant connection $\Omega_W$(Definition 2.3) ,which can
be used to deform every finite dimensional representation of $Br_W
(\Pi )$ to a one parameter class of representations of the
associated braid group $B_W$. In this paper we try to introduce BMW
algebras for every Coxeter group, as deformations of the Brauer type
algebras of \cite{CH}. The starting point of our construction is the
following observation.

 Let $\Gamma$ be $A_2$, the Dynkin diagram with two vertex $v_0 ,v_1
 $ and a simple lace connecting them. So the BMW algebra $B_{A_2}(\tau
 ,l)=B_3 (\tau ,l)$ is generated by $\{ X_0 ,X_1 , E_0 , E_1  \}$ with relations
 as in above table. Denote the free monoid generated by  $\{ X_0 ,X_1 , E_0 , E_1
 \}$
 as $\Lambda$.

 Then for $i=0,1$, there is a function $\Phi ^i $
 on $\Lambda$ such that the following relation holds in  $B_{A_2}(\tau
 ,l)$ :
\begin{equation}
  E_i x  E_i = \Phi ^i (x) E_i
\end{equation}

   These relations are clear from the diagram presentation of the BMW algebras.
It is an interesting fact that we can obtain above $\Phi ^i$ from
  the LK representation of the braid group $B_3$. Let $\sigma_0
  , \sigma_1$ be the generators of $B_3$ in its canonical
  presentation . The LK representation is a 3 dimensional
  representation with 2 parameters $\{ \tau ,l \}$ $\rho_{LK}: B_3 \rightarrow GL(V)$. We can show (Lemma 3.1)
  that for $i=0,1$, the operator $e_i = \frac{l}{\nu} (\rho (\sigma_i )^2 + \nu \rho (\sigma_i
  )-1)$  is a  projector to some line $\mathbb{C} u_i $ in $V$, where $\nu=\frac{l-l^{-1}}{1-\tau }$.   The
  LK representation can be looked as a representation of $B_{A_2} (\tau ,l)
  $ by sending $X_i , E_i$ to $\rho_i (\sigma_i ) , e_i $
  respectively for $i=0,1$. We denote this representation by the same symbol $\rho_{LK}$.
  For any word $W= A_1 A_2 \cdots A_k $ in $\Lambda$,where $A_i \in \{ X_0 ,X_1 , E_0 , E_1
 \}$ we set $\rho_{LK} (W) = \rho_{LK}(A_1) \cdots \rho_{LK}(A_k ) \in
 gl(V)$.  Then there are functions $\bar{\Phi  }^i$ ($i=0,1$) on
 $\Lambda$ such that $e_i \rho_{LK} (W) e_i = \bar{\Phi }^i (W)
 e_i$, since $e_i$ is a projector. Comparing above two identities, we have
$\Phi ^i =\bar{\Phi  }^i $.

Denote the Dihedral group of type  $I_2 (m)$ as $D_m$, whose Dynkin
diagram consists of two vertexes $v_0 , v_1$ and an edge with weight
$m$ connecting them. Denote the associated Artin group as $A_{D_m}$.
The group $A_{D_m}$ has a canonical presentation with two generators
$\sigma_0 ,\sigma_1$.

The group $D_m$ has a natural representation on $\mathbb{C}^2$ as a
reflection group of order $2m$. It has $m$ reflection hyperplanes
$H_0 , H_1 , \cdots , H_{m-1} $. Denote the reflection by $H_i$ as
$s_i$, and chose a linear function $f_i$ on $\mathbb{C}^2$ such that
$H_i =ker f_i$. Set $\omega_i = d f_i / f_i$, which is a holomorphic
1-form on $M$. Denote the complementary space $\mathbb{C}^2
\setminus \cup_{i=0} ^{m-1} H_i $ as $M_m $. The Artin group
$A_{D_m}$ also has a generalize Lawrence-Krammer representation as
special cases of the one constructed by Marin \cite{Ma2} as follows.
Set $V_m = \mathbb{C}<v_0 , v_1 ,\cdots ,v_m
> $ be a vector space with a basis in one to one correspondence with
the set of reflection hyperplanes of $D_m$. The action of $D_m$ on
the second set naturally induces a representation of $D_m$ on $V_m$
which will be denoted as $\iota$. Then on the trivial bundle $M_m
\times V_m$, we have Marin's  connection: $\Omega_{LK} = \Sigma
_{i=0} ^{m-1} \kappa (k \iota (s_i) -p_i ) \omega_i $, where $p_i
\in End(V_m )$ is suitable projector to the line $\mathbb{C} v_i $.
It is proved \cite{Ma2} that $\Omega_{LK}$ is flat and $D_m -$
invariant. So it induces a flat connection on the flat bundle
$M\times _{D_m} V$, whose monodromy gives the generalized LK
representation of $A_{D_m}$, which will be denoted as $T_{LK} :
A_{D_m} \rightarrow Aut (V_m )$.

First we suppose the cases when $m$ is odd, we
have\\

\nd{\bf Lemma 3.1} {\it Set $q=\exp(\kappa k\pi \sqrt{-1})$ and $l=
\exp (\kappa \alpha \pi \sqrt{-1}) / q$. For generic $\kappa$,

$(1)$ $(T_{LK}(\sigma_i )- l^{-1}  ) (T_{LK}(\sigma_i) - q
)(T_{LK}(\sigma_i )+ q^{-1} ) =0$ for $i=0,1$.

$(2)$  $e_i = \frac{l}{q^{-1}-q} (T_{LK}(\sigma_i) - q)
(T_{LK}(\sigma_i) + q^{-1})$ is a projector, i.e, $Rank e_i =1$ for
$i=0,1.$ }\\

Still let $\Lambda$ be the monoid freely generated by $X_0 ,X_1 ,
E_0 ,E_1$. Define a morphism $f: \Lambda \rightarrow End(V_m )$ by
setting:  $f(X_i )= T_{LK}(\sigma_i )$, $f(E_i )= e_i $ $(i=0,1)$.
Considering $(2)$ of above Lemma 3.1 , we define a function
 $ \Phi^i _{m,\Pi ,\kappa}: \Lambda \rightarrow \mathbb{C} $ by:  $e_i f(X) e_i = \Phi^i _{m,\Pi}(X) e_i
 $. Then we define the BMW algebra $B_{D_m } (\Pi ,\kappa)$ associated with $D_m$
 as the following Definition 3.4.

 Suppose $k\in \mathbb{C}\setminus \{ 0 \}$, $\alpha \in \mathbb{C}$, denote $\Pi=\{ k, \alpha \}$. Set
$q=\exp(\kappa \sqrt{-1}\pi k )$ and $l=\exp(\kappa \alpha \pi \sqrt{-1})/q$.\\

\nd{\bf Definition 3.4} {\it The algebra $B_{D_m}(\Pi ,\kappa)$ is
generated by $X_0 ,X_1 ,E_0 ,E_1$ with following relations.

$1)$ $X_i X_i ^{-}=1$, $i=0,1$;

$2)$ $[X_0 X_1 \cdots ]_{m}=[X_1 X_0 \cdots ]_{m} $;

$3)$ $(q ^{-1} -q) E_i = l (X_i -q )(X_i + q ^{-1} )$ for $i=0,1$;

$4)$ $X_i E_i = E_i X_i = l ^{-1} E_i $ for $i=0,1$;

$5)$ $E_i X E_i = \Phi^i _{m, \Pi ,\kappa} (X) E_i$ for any $X\in
\Lambda$,$i=0,1$. }\\

The relation $1) $ to $4)$ have clear analogues with the BMW
 algebra $B_3 (\tau ,l)$. Relation $5)$ is similar with the equation (1) for $B_3 (\tau ,l)$. From above definition we see $B_{D_m} (\Pi ,\kappa)$ is a
 quotient of the group algebra $\mathbb{C} A_{D_m } $, and has the
 Hecke algebra of type $W$ as its quotient.

  When $m$ is even, we have a pair of projectors $e_0 ,e_1$ acting on
  Marin's generalized LK representation as in Lemma 3.2. By using
  them we define a function $\Phi^i _{m,\Pi,\kappa} : \Lambda \rightarrow
  \mathbb{C}$, which is used in Definition 3.3 of the BMW algebra $B_{D_m}
  (\Pi ,\kappa
  )$.

  In section 4 we study the algebras $B_{D_m} (\Pi ,\kappa )$, the following
  results support that they are indeed the deformation of the Brauer
  type algebras $Br_{D_m} (\Pi )$.  Our Theorem 4.1 says that for
  generic parameters $\Pi $, $\dim Br_{D_m} (\Pi) = \dim B_{D_m }(\Pi
  ,\kappa)$.  Proposition says the representations of  $A_{D_m}$
 obtained from monodromy of the KZ connection $\Omega_{D_m}$ (associated with finite
 dimensional representations
of $Br_{D_m}(\Pi )$ factor through $B_{D_m} (\Pi, \kappa)$.

Once the definition of BMW algebras associated with Dihedral groups
are given, they induce a natural definition of BMW algebra
$B_{W_{\Gamma}} (\Pi ,\kappa )$ for any Coxeter group $W_{\Gamma}$,
as in Definition 4.1.

We hope these algebras $B_{W_{\Gamma}}(\Pi ,\kappa)$ be deformation
of the Brauer type algebras $Br_{W_{\Gamma}} (\Pi )$ as the dihedral
cases. We can't give a proof of them in this paper. We only prove
(Theorem 5.1) that $B_{W_{\Gamma}}(\Pi ,\kappa)$ is finite
dimensional when $W$ is a finite group, by similar arguments used in
\cite{CGW1} in the proof of Theorem 1.1 of that paper.  Proposition
5.2 says the algebra $B_{W_{\Gamma}} (\Pi ,\kappa  )$ is isomorphic
to the simply laced BMW algebra defined in  \cite{CGW1} when $W$ is
a simply laced Coxeter group. At last, Theorem 5.2 says for finite
$W$, for generic parameters $\Pi$, the representations of $A_W$
resulted from monodromy of the KZ connection $\Omega_W$ all factor
through $B_W (\Pi ,\kappa )$ for suitable parameters. \\
\nd {\bf Acknowledgements} I 'd like to present my thanks to Ivan
Marin for his detailed explanation of some tools involved in proof
of Lemma 3.5.

\section{preliminaries}

Let $E$ be a complex linear space. An hyperplane arrangement (or
arrangement simply ) in $E$ means a finite set of hyperplanes
contained in $E$. Let $\mathscr{A} =\{ H_i \} _{i\in I}$ be an
arrangement in $E$, we denote the complementary space $ E- \cup
_{i\in I} H_i $ as $M_{\mathscr{A}}$.  Intersection of any subset of
$\mathscr{A}$ is called an edge. If $L$ is an edge of $\mathscr{A}$,
define $\mathscr{A} _L = \{ H_i \in \mathscr{A} | L\subset H_i \}
=\{H_i \} \ and\ I_L = \{ i\in I | L\subset H_i  \} .$

For every $i\in I$, chose a linear form $f _i $ with kernel $H_i $.
Set $\omega _i = d\log f_i $, which is a holomorphic closed 1-form
on $M_{\mathscr{A}}$.  Consider the formal connection $\Omega =
\kappa \sum _{i\in I} X_i \omega _i $.  Here $X_i $ are linear
operators to be determined. When we take $X_i $'s as endomorphisms
of some linear space $E$, then $\Omega $ is realized as a connection
on the bundle $M_{\mathscr{A}} \times E$.  We have the following
theorem of Kohno.

\begin{thm}[Kohno \cite{Ko1}]
 The formal connection $\Omega $ is flat if and only if:

 $[ X_i ,\sum _{j\in I_L } X_j ] =0$ for any codimension  2  edge  $L$  of  $\mathscr{A}$,and for any
$i\in I_L$. Where $[A,B]$ means $AB-BA$.
\end{thm}

Let $W\subset U(E)$ be a finite pseudo reflection group. Let $R$ be
the set of pseudo reflections contained in $W$£¬ and let $\mathcal
{A}=\{ H_i \}_{i\in P}$ be the set of reflection hyperplanes in $E$.
Since the action of $W$ on $E$ send reflection hyperplanes to
reflection hyperplanes, we have a natural action of $W$ on $P$ which
will denoted simply by $(w, i)\mapsto w(i)$ for $w\in W, i\in P$.
For $s,s^{'} \in R$, denote $s\sim s^{'}$ if they lie in the same
conjugacy class. For $i,i^{'} \in P$, denote $i\sim i^{'}$ if there
exist $w\in W$ such that $w(i)=i^{'}$. For $i,j\in P$, denote
$R(i,j)= \{ s\in R | s(H_j)= H_i \}$.  For $i\in P$, let $W_i$ be
the subgroup of $W$ consisting of elements that fixing $H_i $
pointwise. Let $m_i = |W_i|$. It is well known that $W_i$ is a
cyclic group and we denote the element in $W_i$ with exceptional
eigenvalue $e^{\frac{2\pi \sqrt{-1}}{m_i}}$ as $s_i$. For $s\in R$,
define $i(s)\in P$ by requesting $H_{i(s)}$ to be the reflection
hyperplane of $s$. Set $M_W =E-\cup_{i\in P} H_i $. The action of
$W$ on $M_W$ is free, by a famous theorem of Steinberg.  The groups
$B_W = \pi_1 (M_W /W)$ and $P_W = \pi_1 (M_W )$ are called the
complex braid group and the complex pure braid group associated with
$W$ respectively.  For $i\in P$, suppose $\omega_i$ is the 1-form on
$M_W$ associated with $H_i$ as above. Suppose $(V,\rho)$ is a
complex linear representation of $W$. Suppose $\Omega = \sum_{i\in
P} X_i \omega_i$ is a connection on the trivial vector bundle $M_W
\times V$, where $X_i \in End(V)$ for any $i\in P$.  In this case we
can get a vector bundle on $M_W /W$ as the orbit space of the free
action of $W$ on $M_W \times V$ : $w \cdot (p ,v )= (w\cdot p , \rho
(w) (v ) )$, for any $w\in W$, $p\in M_W$ and $v\in V$. Denote the
quotient bundle as $M_W \times _{W} V$, and the quotient map from
$M_W \times V$ to $M_W \times _{W} V$ as $\pi$. We have the
following Lemma, for details and backgrounds about these connections
see Section 4 of \cite{BMR} .

\begin{lem} If $\rho(w) X_i \rho(w)^{-1} = X_{w(i)}$ for any $i\in P$ and $w\in W$, then $\Omega$ induce a connection $\bar{\Omega}$ on
$M_W \times _{W} V$, whose pull back through by $\pi$ is $\Omega$.
\end{lem}

If a connection satisfy the condition in above Lemma, we call it as
$W-$invariant.

Suppose $E^{\mathbb{R}}$ is a $n$ dimensional Eclidean space, denote
the bilinear form on $E^{\mathbb{R}}$ as $<,>$.   $W\subset
O(E^{\mathbb{R}})$ is a finite reflection group with the set of
reflections denoted as $R$. Suppose $W$ is essential, which means
the fixed subspace of $W$ in $E^{\mathbb{R}}$ is $0$. Denote the
complexification $E^{\mathbb{R}} \oplus \sqrt{-1}E^{\mathbb{R}}$ of
$E^{\mathbb{R}}$ as $E$. There is a natural positive definite
Hermitian form on $E$ extending the metric on $E^{\mathbb{R}}$
defined as:  $$< v+\sqrt{-1}w,v^{'}+\sqrt{-1}w^{'}  >= <v,v^{'}> +
\sqrt{-1} <w, v^{'}> -\sqrt{-1}<v,w^{'}> + <w,w^{'}>.$$  Denote the
distance of two points $p,q \in E$ determined by this Hermitian
metric as $d(p,q)$. Because $R$ contain only order 2 elements, the
map $i$ from $R$ to $P$ (index set of reflection hyperplanes) is
bijective. So we simply set $P=R$ and denote the reflection
hyperplane in $E^{\mathbb{R}}$ of $s$ as $H_s$. Denote the
complexification of $H_s$ as $H^{\mathbb{C}} _s$, which is a
hyperplane in $E$.  It is a well-known fact (for example see
\cite{Hu}) that each path component of $E^{\mathbb{R}}-\cup_{s\in R}
H_s $ is a simplicial cone. Choose one from these cones, denote the
closure of it as $\Delta$.  Denote the $n$ reflection hyperplanes
(of $W$) corresponding to $n$ codimension one faces of $\Delta$ as
$H_{s_1}, \cdots ,H_{s_n }$. The following Theorem is canonical (see
\cite{Hu}).

\begin{thm}
The group $W$ is generated by $s_1 , \cdots ,s_n$ with the following
relations.

1) $s^2 _i =1$, $1\leq i\leq n$;  $\ \ \ \ $2) $[s_i s_j \cdots
]_{m_{i,j}} = [s_j s_i \cdots ]_{m_{i,j}}$, $1\leq i\neq j \leq n$.
\end{thm}

Here the expression $[xy\cdots]_{m}$ means a word with length $m$
,in which the letters $x,y$ appear alternatively, and whose left
most subword is $xy$. Similarly we define $[\cdots xy]_{m}$.

The symmetric matrix $\Gamma = (m_{i,j})_{n\times n}$ is called the
Coxeter matrix of $W$, which is equivalent to a Dynkin diagram under
certain simple rules.

Now consider the space $M_W /W$.  Denote the quotient map from $M_W$
to $M_W /W$ as $\pi$. Choose a inner point $p\in \Delta $. Notice
since $E^{\mathbb{R}}$ is a subspace of its complexification $E$,
$p$ belongs to $M_W$. Denote $\pi(p)$ as $\bar{p}$.  Denote the
segment connecting $p$ with $s_i (p)$ as $\gamma_i $, which
intersecting $H_{s_i}$ in one point $p_i$.  Chose $0< \epsilon<<1$.
Let $p^{'}_i$ be the point on $\gamma_i \cap \Delta$ satisfying:
$d(p^{'}_i ,p_i )=\epsilon$. Set $p^{"}_i = s_i (p^{'}_i)$ (also
belong to $\Gamma_i $). Denote the vector $s_i (p) -p$ as
$\alpha_i$, then the reflection $s_i$ can be presented as $s_i (v)=
v- \frac{2<v,\alpha_i
>}{<\alpha_i ,\alpha_i
>} \alpha_i $. denote the path from $p$ to $p^{'} _i$ as
$\gamma^{'}_i$, and the path from $p^{"}_i$ to $s_i (p)$ as
$\gamma^{"}_i$. Define the following path $\eta_i$ from $p^{'} _i$
to $p^{"}_i$ as: $\eta_i : [0,1]\rightarrow E$, $t\mapsto p^{'} _i -
(1-\exp {\pi \sqrt{-1}t}) \frac{<p^{'} _i , \alpha_i  >}{<\alpha_i
,\alpha_i >}\alpha_i$. Denote the composed path $\gamma^{"}_i \circ
\eta_i
 \circ \gamma^{'}_i$ as $l_i$. It is easy to see when $\epsilon$ is
 small enough, the path $l_i$ lie in $M_W$, and it project to a
 closed path  $\bar{l}_i$ in $M_W /W$ based at $\bar{p}$. Denote the
 element of $\pi_1 (M_W /W ,\bar{p} )$ represented by $\bar{l}_i$ as
 $\sigma_i $.

 \begin{thm}[Brieskorn \cite{Bri}]
The group $\pi_1 (M_W /W ,\bar{p} )$ is generated by $\sigma_1
,\cdots ,\sigma_n $ with the following relations.

 $[\sigma_i \sigma_j \cdots ]_{m_{i,j}}= [\sigma_j \sigma_i \cdots
 ]_{m_{i,j}}$ for $1\leq i\neq j\leq n$.

 \end{thm}

 The group $\pi_1 (M_W /W ,\bar{p}  )$ is the Artin group of type $W$, which will be denoted as $A_W$.
  The following results about monodromy of flat connections will be
 used later. Suppose $W$ is a finite Coxeter group, and $M_W , \Delta , s_i , H_i
 ,p,\bar{p}$ defined as above.  Suppose $\rho : W\rightarrow GL(V)$(or $Aut(V)$) is a finite dimensional complex linear
 representation of $W$, and $\Omega =\kappa \sum_{s\in R} X_s \omega_s
 $ is a $W$-invariant flat connection on a trivial vector bundle $M_W \times
 V$, where $X_s \in End(V)$ for any $s\in R$. Denote the flat
 connection on $M_W \times_W V$ induced by $\Omega$ as
 $\bar{\Omega}$. Denote the monodromy representation of
 $\bar{\Omega}$ as $T: \pi_1 (M_W /W , \bar{p}) \rightarrow GL(V)$.
 Then

 \begin{lem}[Proposition 2.3 of \cite{Ma2}] For generic $\kappa$,
 the operator $T(\sigma_i )$ is conjugated to $\rho(s_i ) \exp (\kappa \sqrt{-1}\pi X_{s_i }
 )$ in $GL(V)$.
 \end{lem}

 Now suppose $L$ is an edge of the arrangement $\{ H_s \}_{s\in R}$,
 which means $L\subset E$ is the nonempty intersection of several members of  $\{ H_s \}_{s\in
 R}$. Set $R_L = \{ s\in R | L\subset H_s \} $, and $M_{W_L } = E -\cup _{s\in R_L } H_s
 $. It is a canonical result (see \cite{Hu}) that the subgroup $W_L$  generated by $R_L $
  is the maximal subgroup of $W$ fixing $L$
 pointwise.  As a reflection group it has its own associated Artin
 group $A_{W_L }$ and pure Artin group $P_{W_L}$.  It is clear that we have identifications:
 $P_{W_L}\cong \pi_1 (M_{W_L} )$ and $A_{W_L} \cong \pi_1 (M_{W_L}
 /W_L)$.  Now choose a  point $q\in L$ such that $q\notin H_s $ for
 any $s\notin R_L$. We can choose $\epsilon>0$ small enough such
 that  the ball $D_{\epsilon} (q) = \{ a\in E | d(a,q)<\epsilon  \}$
 has empty intersection with any $H_s (s\notin R_L )$.  It is easy
 to  see the natural morphism $\pi_1 (D_{\epsilon} (q) \cap M_W ) \rightarrow \pi_1
 (M_{W_L})$ induced by the inclusion map is an isomorphism, thus the
 morphism $\pi_1 (D_{\epsilon} (q) \cap M_W) \rightarrow \pi_1 (M_W ) $
 induced by inclusion map gives us an injective morphism: $\lambda_L: P_{W_L}\rightarrow
 P_{W}$. Since $D_{\epsilon} (q)$ is setwise stabilized by $W_L$,
 $\lambda_L$ can be extend to an injective morphism $\bar{\lambda}_L: A_{W_L}\rightarrow
 A_W$.

 Assume $(V,\rho ),\  \Omega$ as above.   On the trivial bundle $M_{W_L} \times V$ we define the following
 connection: $\Omega^{'} =\kappa \sum_{s\in R_L} X_s \omega_s $.
 it is easy to see $\Omega^{'}$ is also a flat and $W_L$-invariant
 connection, so it induce a flat connection $\bar{\Omega}^{'}$ on
 the quotient bundle $M_{W_L} \times_{W_L} V$. The connection
 $\bar{\Omega}^{'}$ induce a monodromy representation $T_L : A_{W_L} \rightarrow
 GL(V)$. On the other hand, by the injection $\bar{\lambda}_L$ we
 have another representation $T\circ \bar{\lambda}_L$ of $A_{W_L}$:
 $A_{W_L}\rightarrow A_{W} \rightarrow GL(V)$. In \cite{Ma2}
  the following theorem is presented.

  \begin{thm}
Assume above conventions. For generic $\kappa$, the representation
$(V , T_L )$ is isomorphic to $(V,T\circ \bar{\lambda}_L )$.
\end{thm}

Let $W\subset U(E)$ be a finite pseudo reflection group. In
\cite{CH} the author introduced a Brauer type algebra $Br_W (\Pi )$
associated with $W$, which is defined as follows.  Set notations
$R,\mathcal {A}, P, i(s), s\sim s^{'}, i\sim i^{'}, R(i,j)$ as
above.  Choose constants $\{ k_s  \}_{s\in R} \cup \{ \alpha_i
\}_{i\in P}$ such that $k_s = k_{s^{'}}$ if $s\sim s^{'}$ and
$\alpha_i =\alpha_{i^{'}}$ if $i\sim i^{'}$. Denote $\{ k_s \}_{s\in
R} \cup \{ \alpha_i \}_{i\in P}$ by one symbol $\Pi$.  For $i\neq j$
we denote $H_i \pitchfork H_j$ if $\{ k \in P \ |\ H_i \cap H_j
\subset H_k \} =\{ i,j \}$.  A codimension 2 edge $L$ of $\mathcal
{A}$ will be called a crossing edge if there exists $i,j\in P$ such
that $H_i \pitchfork H_j$ and $L= H_i \cap H_j$, otherwise $L$ will
be called a noncrossing edge.


\begin{defi}
  The algebra $Br_W (\Pi )$
associated with pseudo reflection group $W\subset U(E)$ is generated
by the set
  $\{ T_{w} \}_{w\in G } \cup  \{ e_i \}_{i\in P}$ which satisfies the following
relations.

$0)$ $T_{w_1} T_{w_2} = T_{w_3 }$ if $w_1 w_2 =w_3 $.

$1)$ $ T_{s_i } e_i  = e_i T_{s_i } = e_i $, for $i\in P$.

$1)^{'}$ $T_w e_i =e_i T_w = e_i $, for $w\in G$ such that $w(H_i) =
H_i$, and  $H_i \cap E_w $ is a noncrossing edge. Where $E_w =\{
v\in E | w(v)=v  \}$.

$2)$ $e_i ^2 = \alpha_i e_i $ .

$3)$ $T_w e_j  = e_i T_w $ , if $w \in G  $ satisfies $w(H_j )=H_i$.

$4)$ $e_i e_j = e_j e_i $, if  $H_i \pitchfork H_j$ .

$5)$ $ e_i e_j = (\sum _{s\in R(i,j) } k_s T_s )e_j = e_i (\sum
_{s\in R(i,j) } k_s T_s )$ ,
 if $H_i \cap H_j$ is a noncrossing edge,
  and $R(i,j) \neq \emptyset $.

$6)$ $e_i e_j =0 $, if $H_i \cap H_j$ is a noncrossing edge, and
$R(i,j) =\emptyset $.

\end{defi}

It is proved in \cite{CH}(Theorem 8.4) that when $W$ is a finite
type Coxeter group with Coxeter matrix $(m_{i,j})_{n\times n}$, the
algebra $Br_W (\Pi )$ is isomorphic to the following algebra
$Br^{'}_{W_M} (\Pi)$  with a canonical presentation.

\begin{defi}
For any Coxeter matrix $M= (m_{i,j})_{n\times n}$, the algebra
$Br^{'}_{W_M} (\Pi)$ is defined as follows. Let $\Pi$ be as in
Definition 2.1. Denote $\alpha_{s_i}$ in $\Upsilon$ as $\alpha_i$.
If we don't give range for an index then it means "for all". The
generators are $S_1 ,\cdots ,S_n ,E_1 ,\cdots ,E_n $. The relations
are
\begin{align*}
&1) S_i ^2 =1 ; &8)&  E_i w E_j =0  \mbox{ for any word } w  \\
& 2) [S_i S_j \cdots ]_{m_{i,j }} =[S_j S_i \cdots ]_{m_{i,j}}; & &
\mbox{composed from } \{S_i ,S_j
\} \mbox{If } m_{i,j}= 2k>2 ; \\
& 3)  S_i E_i =E_i =E_i S_i ;  & 9)&  E_i [S_j S_i \cdots ]_{2l-1 }
E_i = (k_s +k_{s^{'}} )
E_i  \\
& 4)  E_i ^2 = \alpha_i E_i ;  & &\mbox{for }  1\leq l\leq k, \mbox{If } m_{i,j}= 2k>2 .  \\
&  5)   S_i E_j =E_j S_i  \mbox{if } m_{i,j } = 2 ; & &\mbox{Where }
s=[S_j S_i \cdots ]_{2l-1}, s^{'}=[S_j S_i \cdots ]_{2(k+l)-1}.\\
& 6)   E_i E_j =E_j E_i  \mbox{if } m_{i,j } = 2 ; &10)& E_i [S_j S_i \cdots ]_{2l-1 } E_i =k_{s_{\epsilon }} E_i  \\
& 7)  [S_j S_i \cdots ]_{2k -1} E_i & &\mbox{for } 1\leq l\leq k, \mbox{if } m_{i,j} = 2k+1 ;  \\
& = E_i [S_j S_i \cdots ]_{2k
-1} =E_i , & &\mbox{Where } \epsilon=i(j) \mbox{if } l \mbox{is odd (even) }.  \\
& \mbox{if } m_{i,j}=2k>2 ;    &11)& [S_i S_j \cdots ]_{2k } E_i =
E_j [S_i S_j \cdots ]_{2k
} \mbox{If } m_{i,j} = 2k+1 .  \\
\end{align*}

\end{defi}
An important feature of the algebra $Br_{W} (\Pi)$ is it supports
the following $W$-invariant, flat formal connection
\begin{defi}The KZ connection of $Br_{W} (\Pi)$ is the following formal connection on
$M_W \times Br_W (\Pi)$: $ \Omega_W =\kappa \sum_{i\in P} ( \sum_{s:
i(s)=i} k_s T_s- e_i   ) \omega_i .  $
\end{defi}
The following fact is from  Proposition 5.1 of \cite{CH}.
\begin{prop} The connection $\Omega_W $ is flat and $W$-invariant.
\end{prop}
 This connection can deform every
finite dimensional representation of $Br_{W} (\Pi) $ to a
representation of $B_W$,the associated braid group. We have the
following example. Let $R, P, M_W , \omega_i $ be defined as above.
Set $V_W = \mathbb{C} < v_i
>_{i\in P}$ be a complex linear space with a basis in one to one
correspondence with the set of reflection hyperplanes of $W$. The
permutating action of $W$ on $P$ induces a representation $\iota $
of $W$ on $V_W$ in a natural way: $\iota(w) (v_i) = v_{w(i)}$ for
any $w\in W$ and $i\in P$. Let the data $\{ k_s  \}_{s\in R} \cup \{
\alpha_i \}_{i\in P}$ be $\Pi$ in $Br_W (\Pi)$.  For $i\in P$,
define an element $p_i \in End(V_W )$ as by $p_i (v_i) = \alpha_i
v_i $; $p_i (v_{j})= (\sum_{t\in R: t(j) =i } k_{t} ) v_i$. So $p_i$
is a projector to the line $\mathbb{C}v_i$.
\begin{lem} The map $T_w \mapsto
 \iota(w) $ for $w\in W$ and $ e_i \mapsto p_i$ for $i\in P$ extends
 to a representation $\rho_{LK}: Br_W (\Pi) \rightarrow End(V_W ) $.
\end{lem}
 We call this representation as the infinitesimal Lawrence-Krammer representation.  So the  connection $\Omega^{LK} _{W} =\kappa \sum_{i\in P}  (
\sum_{s\in R: i(s)=i} k_s \iota(s) - p_i ) \omega_i $ on the vector
bundle $M_W \times V_W$ is flat and $W$-invariant by Proposition
2.1.   $\Omega^{LK} _{W}$ induces a flat connection
$\bar{\Omega}^{LK} _W$ on the bundle $M_W \times _{W} V_W $ ,whose
monodromy representation is defined as the generalized
Lawrence-Krammer representations (of $A_W $).  The generalized
Lawrence-Krammer representations for  finite type simply-laced Artin
groups are invented  by Cohen and Wales in \cite{CW}, and Digne in
\cite{Di} independently. In \cite{Ma2} Marin introduced  generalized
Lawrence-Krammer representations for pseudo reflection groups whose
pseudo reflections all have order 2.  Above slightly further
generalization of Marin's work can be found in \cite{CH}, which is
put in the frame of $Br_W (\Pi)$ for later convenience.

Consider the permutating action of $W$ on $P$, denote the space of
orbits as $P/W$. For each $c\in P/W$, it is easy to see
$V_{W,c}=\mathbb{C}< v_i >_{i\in c}$ is a sub representation of
$\rho_{LK}$. Thus we have a decomposition $\rho_{LK} = \oplus_{c\in
P/W} \rho_{LK,c}$ according to the decomposition $V_W = \oplus_{c\in
P/W} V_{W,c}$. In dihedral cases (where we can set $P=R$ ), when $m$
is odd, the set of reflections $R$ of $D_m$ contains only one orbit.
When $m$ is even, let $c_0 \subset R$, $c_1 \subset
 R$ be those reflections conjugated to $s_0$, $s_1$ respectively.
 Then we have $V_{D_m} = V_{D_m , c_0} \oplus V_{D_m ,c_1}$,
 correspondingly $\rho_{LK} = \rho_{LK,c_0} \oplus \rho_{LK,c_1}$.

 The following result follows from
Proposition 3.4 and Proposition 4.2 of \cite{Ma2}.
\begin{lem} When $m$ is odd, for generic $k,\alpha$, the
infinitesimal Lawrence-Krammer representation $\rho_{LK}$ is
irreducible. When $m$ is even,  for generic $k_0 ,\alpha_0 ,k_1
,\alpha_1$, the representations $\rho_{LK,c_0}, \rho_{LK,c_1}$ are
irreducible.
\end{lem}

 \section{Generalized Lawrence-Krammer Representations}
Suppose $W$ is a finite Coxeter group. We assume the conventions for
$E, R, M_W  $ as in Section 2. We study the generalized
Lawrence-Krammer representation in detail for the cases when $W$ is
a dihedral group, as they will play an important role in the next
section.

Denote the dihedral group of type $I_2 (m)$ as $D_m $. The
arrangement of its reflection hyperplanes can be explained with the
following Figure 1.
\begin{figure}[htbp]

  \centering
  \includegraphics[height=4.5cm]{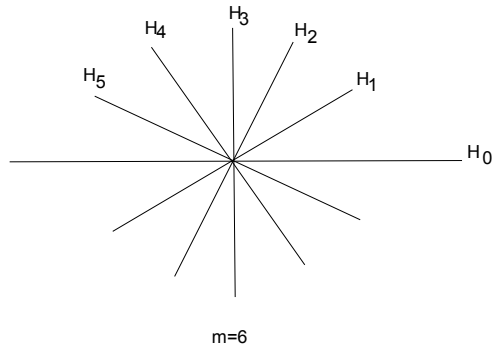}
  \caption{The arrangement of $D_6$ }
\end{figure}

There are $m$ lines(hyperplanes) passing the origin. The angle
between every two neighboring lines is $\pi /m$. Suppose the
$x$-axis is one of the reflection lines and denote it as $H_0 $, we
denote these lines by $H_0 , H_1 ,\cdots , H_{m-1} $ in
anticlockwise order as shown in above graph.  Denote the reflection
by $H_i $ as $s_i $. The set of reflections in $D_m$ is $R=\{ s_i \}
_{0\leq i\leq m-1 }$. Denote the rotation of $2j \pi /m $ in
anticlockwise order as $r_j $. It is well known that $D_m$ is
generated by $s_0 ,s_1 $ with the following presentation
$$< s_0\ ,s_1 , \ | ( s_0 s_1 )^m =1 , s_0 ^2 =s_1 ^2 =1 >.$$
The associated Artin group $A_{W_m}$ has a presentation:
$$<\sigma_0 ,\ \sigma_1 ,\ | [\sigma_0 \sigma_1 \cdots ]_m =[\sigma_1 \sigma_0 \cdots ]_m >.$$
As a set $D_m =\{ s_i ,r_i \}_{0\leq i\leq m-1}$.  Under this
presentation, $s_i $ can be determined inductively in the following
way.  $s_1 =s_1 $, $s_2 = s_1 s_0 s_1 $ and $s_i = s_{i-1} s_{i-2}
s_{i-1} $. Recall by $[s_i s_j \cdots ]_{k}$ we denote the length
$k$ word starting with $s_i s_j $, in which $s_i $ and $s_j $ appear
alternatively. The word $[\cdots s_0 s_1]_k$ is defined similarly.
Then $s_i = [s_1 s_0 \cdots ]_{2i-1} $, where $s_m =[s_1 s_0 \cdots
]_{2m-1 } =s_0 $.\\

\nd{\bf  The cases when $m$ is odd } Recall the connection
$\Omega^{LK} _W$ on $M_W \times V_W$ in Lemma 2.2. We denote $D_m$
simply as $W$ in this paragraph. For these cases all reflections in
$W$ lie in one conjugating class. So the constants $k_s$ all equal
to some $k$ and constants $\alpha_s$ all equal to some $\alpha$. We
simply denote the base elements $v_{s_i}$ of $V_W$ as $v_i$,  the
projector $p_{s_i}$ as $p_i$, and $\omega_{s_i}$ as $\omega_i$ for
$0\leq i\leq m-1$. The connection can be rewritten as:
\begin{equation}
\Omega_m = \sum^{m-1} _{i=0}\kappa (k \iota ( s_i ) - p_i )\omega_i
,
\end{equation}
Where the $W$ representation $(V_W ,\iota)$ is as in Section 2, and
$p_i $ is defined by $p_i (v_i )= \alpha v_i ; p_i (v_j)=k v_i $ if
$j\neq i$. Denote the induced flat connection on $M_W \times_W V_W$
as $\bar{\Omega}_m$ and the monodromy representation of
$\bar{\Omega}_m$ as $T_{LK}: A_W \rightarrow GL(V_W)$.
\begin{lem} Set $q=\exp(\kappa k\pi \sqrt{-1})$ and $l= \exp (\kappa \alpha \pi \sqrt{-1}) / q$. For generic $\kappa$,

$(1)$ $(T_{LK}(\sigma_i )- l^{-1}  ) (T_{LK}(\sigma_i) - q
)(T_{LK}(\sigma_i )+ q^{-1} ) =0$ for $i=0,1$.

$(2)$  $e_i = \frac{l}{q^{-1}-q} (T_{LK}(\sigma_i) - q)
(T_{LK}(\sigma_i) + q^{-1})$ is a projector, i.e, $Rank e_i =1$ for
$i=0,1.$
\end{lem}
\begin{proof} We only need to prove the case $i=0$. Another case is
similar.  First we diagonalize $\iota(s_0)$ under certain basis $\{
u_0 ,\cdots, u_{m-1}  \}$ where $u_0 =v_0$. Under the same basis
$\kappa( k \iota(s_0) - p_0 )$ is presented by a matrix whose
entries other than the diagonal and the first column are zero. Now
by Lemma 2.2, the operator $T_{LK}(\sigma_0)$ is conjugated to
$\iota(s_0) \exp(\sqrt{-1}\pi \kappa(k \iota(s_0)-p_0 )) $. Under
the basis $\{ u_i \}$, the operator $\iota(s_0) \exp(\sqrt{-1}\pi
\kappa(k \iota(s_0)-p_0 )) $ is still presented by a matrix whose
entries other than the diagonal and the first column are zero, and
whose $(1,1)$ entry is $l^{-1}$ ,other diagonal entries are $q,
q^{-1}$.  When $\kappa (k- \alpha) \notin \{ \kappa k +\mathbb{Z} \}
\cup \{ -\kappa k + \mathbb{Z} \}
 \}$, so $l^{-1} \notin \{ q, q^{-1} \}$, we see $\iota(s_0) \exp(\sqrt{-1}\pi \kappa (k
\iota(s_0)-p_0 ))$ is also diagonalizable with eigenvalues $q,
-q^{-1} , l^{-1}$, and the eigenspace of $l^{-1}$ is one
dimensional. So complete the proof.

\end{proof}

\nd{\bf  The cases when $m=2k$ is even } In this paragraph denote
$D_m$($m$ being even)as $W$. In these cases the set of reflections
in $W$ consists of two conjugacy classes, represented by $s_0$,
$s_1$ respectively. So constants for the connection consists of
$k_{s_0} ,k_{s_1}, \alpha_{s_0}, \alpha_{s_1}$. Denote them simply
by $k_{0} ,k_{1}, \alpha_{0}, \alpha_{1}$. Denote the basis element
$v_{s_i}$ of $V_W$ as $v_i$, projectors $p_{s_i}$ as $p_i$, and
$\omega_{s_i}$ as $\omega_i$ for $0\leq i\leq m-1$. The connection
$\Omega^{LK} _W $ can be rewritten as
\begin{equation}
\Omega_m =\kappa(  \sum^{k-1} _{i=0} (k_0 \iota(s_{2i} )-
p_{2i})\omega_{2i} + \sum^{k-1} _{i=0} (k_1 \iota(s_{2i+1}) -
p_{2i+1})\omega_{2i+1}).
\end{equation}
 Where $p_i$
$(0\leq i\leq m-1)$ is determined by: $p_i (v_j )= 0$ if $2\nmid
i-j$; $p_i (v_j)=2 v_i $ if $2 \mid i-j$ and $i\neq j$; $p_i (v_i)=
k_0 v_i$ if $i$ is odd; $p_i (v_i) = k_1 v_i$ if $i$ is even. Denote
the flat connection on $M_W \times_{W} V_W$ induced from $\Omega_m$
as $\bar{\Omega}_m$, and denote the monodromy representation of
$\bar{\Omega}_m$ as $T_{LK}: A_W \rightarrow GL(V_W)$. We have the
following Lemma. The proof is similar to Lemma 3.1.
\begin{lem}
Set $q_i =\exp(\kappa k_i \pi \sqrt{-1})$ and $l_i = \exp (\kappa
\alpha_i \pi \sqrt{-1}) /q_i $ for $i=0,1$. For generic $k_i$ and
$\alpha_i$,

$(1)$ $(T_{LK}(\sigma_i )- l_i ^{-1}  ) (T_{LK}(\sigma_i) - q_i
)(T_{LK}(\sigma_i )+ q_i ^{-1} ) =0$ for $i=0,1$.

$(2)$  $e_i = \frac{l_i }{q_i ^{-1}-q_i } (T_{LK}(\sigma_i) - q_i )
(T_{LK}(\sigma_i) + q_i ^{-1})$ is a projector, i.e, $Rank e_i =1$
for $i=0,1.$

\end{lem}

  The following lemma is a easy fact.
\begin{lem} Suppose $V$ is a Hermitian vector space whose distance is denoted as $d(, )$. Suppose $\{ v_i  \}_{i=1,\cdots
,m} \subset V$ span $V$. Then there exists $\epsilon >0$, such that
$d(v^{'} _i ,v_i) <\epsilon $ for $i=1,\cdots,m$ implies $\{ v^{'}_i
\}_{i=1,\cdots ,m}$ span $V$.
\end{lem}
\begin{defi}
Let $V$ be a linear space. A set $\{ A_i   \}_{i\in I}\subset gl(V)$
is called irreducible on $V$ if for any $v\in V$, the set $\{ A_i
(v) \}_{i\in I}$ span $V$.
\end{defi}
\begin{lem}
Let $V$ be a finite dimensional Hermitian linear space. Suppose $\{
 A_i  \}_{i=1, \cdots ,m}\subset gl(V)$ is a finite set being irreducible on $V$.
 Then there exists $\epsilon >0$ such that $d(A^{'} _i , A_i) < \epsilon
 $ for $i=1,\cdots ,m$ implies that $\{ A^{'}_i  \}_{i=1,\cdots ,m}$
 is irreducible on $V$. Where $A^{'} _i \in gl(V)$, and $d(\ ,\ )$
 is a metric on $gl (V)$ defined by $d ( (a_{i,j}), (b_{i,j})  )= \sum _{i,j} |a_{i,j}
 -b_{i,j}|^2$ (Here we suppose a basis of $V$ has been chosen, so identify $gl(V)$ with the matrix algebra.)
\end{lem}
\begin{proof} Denote the projective space of $V$ as $P(V)$. By Lemma
4.3, for any $[v]\in P(V)$ we can easily find a  neighborhood
$U([v])$ of $[v]$ and a constant $\epsilon_{[v]} >0$ such that
$[v^{'}] \in U([v])$ and $d( A^{'} _i ,A_i )<\epsilon $ for
$i=1,\cdots,m$ implies that $\{ A^{'} _i v^{'}  \}_{i=1,\cdots ,m}$
span $V$. From this we can find a constant $\epsilon >0$ such that
if $d( A^{'} _i ,A_i )<\epsilon $ for $i=1,\cdots,m$, then $\{ A^{'}
_i \}_{i=1,\cdots ,m}$ is irreducible on $V$, because $P(V)$ is
compact.

\end{proof}
\begin{lem}When $m$ is odd,  for generic data $\Pi$ and $\kappa$, the representation
$T_{LK}$ is irreducible. When $m$ is even, for generic data $\Pi$
and $\kappa$, the representations $T_{LK,0}, T_{LK,1}$ are
irreducible.
\end{lem}
\begin{proof} First consider the cases when $m$ is odd. As in
section 2, set $\Pi = \{ k,\alpha \}$, the map $s_i \mapsto \iota
(s_i )$, $e_i \mapsto p_i $ $(i=0,1)$ extends to a representation
$\rho_{LK}: Br_{D_m} (\Pi) \rightarrow End(V_{D_m})$. Where $p_i ,
\iota(w)$ are defined as in section 3. Then the connection
$\Omega_m$ (equation (1)) can be rewritten as $\Omega_m = \sum^{m-1}
_{i=0} (k\rho_{LK} (s_i) -\rho_{LK} (e_i))\omega_i$.  By Theorem 6.1
of \cite{CH} we know $Br_{D_m} (\Pi)$ has dimension $2m+m^2$.
Suppose $x_1 ,\cdots , x_{2m+m^2}$ is a basis of $Br_{D_m}(\Pi)$,
and suppose the generating set $\{ s_0 ,s_1 ,e_0 ,e_1 \} \subset
{x_1 ,\cdots , x_{2m+m^2} }$.  By Lemma 2.4, we know for generic
$\Pi$, the set of operators $\{ \rho_{LK} (x_1) ,\cdots , \rho_{LK}
(x_{2m+m^2}) \}$ is irreducible on $V_{D_m}$. Now by Lemma 2.2, we
have \begin{equation} T_{LK} (\sigma_0 ) = \rho_{LK} (s_0) \exp
(\sqrt{-1} \pi \kappa( k \rho_{LK} (s_0) -\rho_{LK} (e_0) ))=
\rho_{LK}(s_0) + O(\kappa) \end{equation} and
\begin{equation}T_{LK} (\sigma_1 ) = A(\kappa) \rho_{LK} (s_1) \exp
(\sqrt{-1} \pi \kappa( k \rho_{LK} (s_1) -\rho_{LK} (e_1) ))
A(\kappa)^{-1}= \rho_{LK}(s_1) +O(\kappa)\end{equation} for suitable
matrix $A(\kappa)$ such that $A(\kappa)= I + O(\kappa)$. So we have
\begin{equation}\frac{1}{2\sqrt{-1}\pi \kappa}( T_{LK}(\sigma_0) ^2 -1) =
k\rho_{LK} (s_0) -\rho_{LK} (e_0) + O(\kappa) . \end{equation}
Similarly we have \begin{equation}\frac{1}{2\sqrt{-1}\pi \kappa}(
T_{LK}(\sigma_1) ^2 -1) = A(\kappa) (k \rho_{LK} (s_1) -\rho_{LK}
(e_1)) A(\kappa)^{-1}  + O(\kappa)  = (k \rho_{LK} (s_1) -\rho_{LK}
(e_1)) + O(\kappa). \end{equation} Apply Lemma 3.4 to the set $\{
\rho_{LK} (x_1) ,\cdots , \rho_{LK} (x_{2m+m^2}) \}$, we obtain some
$\epsilon
>0 $. Since $\{ s_0 ,s_1 ,ks_0 -e_0 , ks_1 -e_1  \}$ generate
$Br_{D_m}(\Pi)$, by above equations $(3), (4), (5), (6)$, we can
find elements $\{ X_1 , \cdots , X_{2m+m^2} \} \subset \mathbb{C}
A_{D_m}$ and $\eta >0$, such that if $\kappa <\eta$ then $d ( T_{LK}
(X_i) , \rho_{LK} (x_i) ) < \epsilon$ for $1\leq i \leq 2m +m^2$.
Where $d(A,B )$ is the distance introduced in Lemma 3.4.  So by
Lemma 3.4 for those $\Pi$ and $\kappa <\eta$, the representation
$T_{LK}$ is irreducible. So we see the subset of $\kappa$ with which
the representations are irreducible contains an open set. Since the
subset of $\kappa$ with which the representation $T_{LK}$ being
reducible is a subvariety, we see for generic $\kappa$ the
representation $T_{LK}$ is irreducible.  The cases when $m$ is even
can be proved similarly.

\end{proof}

Denote the free monoid generated by $X_0  ,X_1 ,E_0 ,E_1$ as
$\Lambda$. Let $m$ be an  odd positve integer, and
$\Pi=(k,\alpha)\in \mathbb{C}^{\times} \times \mathbb{C}$, or $m$
being an even positive integer and $\Pi= (k_0 ,k_1 , \alpha_0 ,
\alpha_1) \in (\mathbb{C}^{\times})^2 \times \mathbb{C}^2 $. For $X
\in \Lambda$, define $f(X)\in End(V_W)$ to be the morphism obtained
by replace letters $X_0 ,X_1 ,E_0 ,E_1$ in $X$ with
$T_{LK}(\sigma_0), T_{LK}(\sigma_1), e_0 ,e_1$ respectively. For
example, $f(X_0 E_0 E_1)= T_{LK}(\sigma_0 ) e_0 e_1 \in End(V_W)$.
\begin{defi}
Choose generic constants system $\Pi$ and $\kappa$.  The map $\Phi^i
_{m, \Pi ,\kappa}: \Lambda \rightarrow \mathbb{C}$ $(i=0,1)$ is
defined by the following identity. $e_i f(X) e_i = \Phi^i _{m, \Pi
,\kappa} (X) e_i$.
\end{defi}
\begin{rem}
Above definition makes sense because for any $\phi \in End(V)$, $e_i
\phi e_i \in \mathbb{C} e_i $ since $e_i$'s are projectors by Lemma
3.1 and Lemma 3.2.
\end{rem}

Now we present the definition of Birman-Murakami-Wenzl algebras for
dihedral groups. Let $W$ be a finite Coxeter group, suppose the set
of reflections $R$ in $W$ consists of $n_W$ conjugacy classes
$\mathcal {C}_i$ $(i=1,\cdots ,n_W )$. We associate a pair of
numbers $(k_i , \alpha_i )\in \mathbb{C}^{\times} \times \mathbb{C}$
with each class $\mathcal {C}_i$, and denote the data of all $(k_i ,
\alpha_i )$ as $\Pi$. We denote the BMW algebra associated with $W$
(to be defined ) as $B_W (\Pi ,\kappa )$, where $\kappa \in
\mathbb{C}$.  If $W$ is of type $\Gamma$ $($ $\Gamma $ is a Coxeter
matrix or Dynkin diagram $)$, we also denote the associated BMW
algebra as $B_{\Gamma}(\Pi ,\kappa )$. Suppose $m=2l$. So $\Pi=\{k_0
,\alpha_0 ,k_1 ,\alpha_1 \} $.  Set $q_i = \exp(\kappa \sqrt{-1}\pi
k_i )$, $l_i = \exp(\kappa \alpha_i \pi \sqrt{-1})/q_i$ for $i=0,1.$
\begin{defi} The algebra $B_{D_m}(\Pi)$ is generated by $X_0 ^{\pm} ,X_1 ^{\pm} ,E_0
,E_1$ with following relations.

$1)$ $X_i X_i ^{-}=1$, $i=0,1;$

$2)$ $[X_0 X_1 \cdots ]_{m}=[X_1 X_0 \cdots ]_{m} $;

$3)$ $(q_i ^{-1} -q_i ) E_i = l_i (X_i -q_i )(X_i + q_i ^{-1} )$ for
$i=0,1$;

$4)$ $X_i E_i = E_i X_i = l_i ^{-1} E_i $ for $i=0,1$;

$5)$ $E_i X E_i = \Phi^i _{m, \Pi ,\kappa} (X) E_i$ for any $X\in
\Lambda$, $i=0,1$;

$6)$ $E_0 X E_1 =0= E_1 X E_0 $ for any $X\in \Lambda$.

\end{defi}

When  is odd, the data $\Pi=\{ k, \alpha \}$. Set $q=\exp(\kappa
\sqrt{-1}\pi k )$ and $l=\exp(\kappa \alpha \pi \sqrt{-1})/q$.
\begin{defi} The algebra $B_{D_m}(\Pi ,\kappa )$ is generated by $X_0 ,X_1 ,E_0
,E_1$ with following relations.

$1)$ $X_i X_i ^{-}=1$, $i=0,1$;

$2)$ $[X_0 X_1 \cdots ]_{m}=[X_1 X_0 \cdots ]_{m} $;

$3)$ $(q ^{-1} -q) E_i = l (X_i -q )(X_i + q ^{-1} )$ for $i=0,1$;

$4)$ $X_i E_i = E_i X_i = l ^{-1} E_i $ for $i=0,1$;

$5)$ $E_i X E_i = \Phi^i _{m, \Pi ,\kappa} (X) E_i$ for any $X\in
\Lambda$,$i=0,1$.

\end{defi}

\begin{rem}
Above $5)$ of Definition 3.3 and $4)$ of Definition 3.4 looks
containing infinite many relations. But it isn't hard to see we can
choose finite many $X$ in these relations to obtain the same
algebra.
\end{rem}

\section{BMW type algebra associated with dihedral groups} In this section we
study the algebras $B_{D_m}(\Pi ,\kappa)$ to certify from several
aspects that they are suitable deformations of the Brauer type
algebras introduced in \cite{CH}. The following Proposition 4.1 is
easy to see.

\begin{prop} Denote the ideal in $B_{D_m}(\Pi ,\kappa )$ generated by $E_0
,E_1$ as $I_m$. Then when $m$ is odd, the quotient algebra
$B_{D_m}(\Pi ,\kappa )/ I_m$ is isomorphic to the Hecke algebra
$H_{D_m} (q)$. When $m$ is even, $B_{D_m}(\Pi ,\kappa )/ I_m$ is
isomorphic to the Hecke algebra $H_{D_m} (q_0 , q_1)$ of
multi-parameters.
\end{prop}
Set $\tau = \frac{l(l^{-1}-q)(l^{-1}+q^{-1})}{q^{-1}-q}; \tau_i
=\frac{l_i (l_i ^{-1}-q_i )(l_i ^{-1}+q_i ^{-1})}{q_i ^{-1}-q_i } $
for $i=0,1$.
\begin{lem}  In $B_{D_m}(\Pi ,\kappa )$ for odd $m$ we have

$(1)$ $E_i ^2 =\tau E_i $ for $i=0,1$;

$(2)$ $E_1 [X_0 X_1 \cdots ]_{m-1} = [X_0 X_1 \cdots ]_{m-1} E_0 $.

In $B_{D_m}(\Pi)$ for even $m$ we have

$(3)$   $E_i ^2 = \tau_i E_i $ for $i=0,1$;

$(4)$  $E_1 [X_0 X_1 \cdots ]_{m-1} = [X_0 X_1 \cdots ]_{m-1} E_1 $;

$(5)$  $E_0 [X_1 X_0 \cdots ]_{m-1} = [X_1 X_0 \cdots ]_{m-1} E_0$;

$(6)$ $[\cdots X_1 X_0]_{m-1} E_1 =E_1 [\cdots X_1 X_0]_{m-1}=
\tau_1 ^{-1} \Phi ^1 _{m, \Pi} ([\cdots X_1 X_0 ]_{m-1}) E_1 ;$

$(7)$ $[\cdots X_0 X_1]_{m-1} E_0 =E_0 [\cdots X_0 X_1]_{m-1}=
\tau_0 ^{-1} \Phi ^0 _{m, \Pi} ([\cdots X_0 X_1 ]_{m-1}) E_0 .$

\end{lem}

\begin{proof}
$(1), (3)$ are by simple computations. To show $(2)$, we use
relation $(3)$ of Definition 3.3 and the identities $X_1 ^i [X_0 X_1
\cdots ]_{m-1} = [X_0 X_1 \cdots ]_{m-1} X_0 ^i $ for nonzero $i$'s.
The proof of $(4),(5)$ are similar to $(2)$. For $(6)$, we start
from the equation $E_1 [\cdots X_1 X_0]_{m-1} E_1 = \Phi ^1 _{m,
\Pi} ([\cdots X_1 X_0 ]_{m-1}) E_1 $, which is a special case of
relation $(5)$ in Definition 3.2. Then by the proved equations $(4),
(1)$ of this lemma, the left side equals $\tau_1 ^{-1} [\cdots X_1
X_0]_{m-1} E_1$ and $(6)$ is proved. The proof of $(7)$ is similar.
\end{proof}

\begin{lem} $(1)$ When $m$ is odd, the algebra $B_{D_m}(\Pi ,\kappa )$ is spanned by the
subset $$\{ [X_0 X_1 \cdots ]_i | 0\leq i\leq m \} \cup \{ [X_1 X_0
\cdots ]_j | 1\leq j\leq m-1 \} \cup \{ [\cdots X_0 X_1 ]_{i} E_0
[X_1 X_0 \cdots ]_j | 0\leq i,j\leq m-1 \}.$$

$(2)$ When $m=2k$ is even, the algebra $B_{D_m}(\Pi ,\kappa )$ is
spanned by the subset $\{ [X_0 X_1 \cdots ]_i | 0\leq i\leq m \}
\cup \{ [X_1 X_0 \cdots ]_j | 1\leq j\leq m-1 \} \cup  \{ [\cdots
X_0 X_1 ]_i E_0 [X_1 X_0 \cdots ]_j | 0\leq i,j\leq k-1 \} \cup \{
[\cdots X_1 X_0 ]_i E_1 [X_0 X_1 \cdots ]_j | 0\leq i,j\leq k-1 \}.
$
\end{lem}

\begin{proof}
First we prove  $(1)$.  It is a well known fact that the Hecke
algebra $H_m (q) $  can be spanned by the subset $\{ [X_0 X_1 \cdots
]_i | 0\leq i\leq m \} \cup \{ [X_1 X_0 \cdots ]_j | 1\leq j\leq m-1
\}$ for any $q$. So by Proposition 4.1, we only need to show the
ideal $I_m$ is spanned by the subset $\{ [\cdots X_0 X_1 ]_{i} E_0
[X_1 X_0 \cdots ]_j | 0\leq i,j\leq m-1 \}$. Now by $(2)$ of Lemma
4.1, $I_m$ can be spanned by $\{ f(X_0 ^{\pm}, X_1 ^{\pm}, E_0   )
\}$, which is the set of words composed by letters $E_0 , X_i
^{\pm}$. By $(5)$ of Definition 3.4 and $(1)$ of Lemma 4.1 we see
$I_m $ can be spanned by those words in $\{ f(X_0 ^{\pm}, X_1
^{\pm}, E_0   ) \}$ which contain only one $E_0$. Then by $(3),(4)$
of Definition 3.3, $I_m$ can be spanned by the set $\{ [\cdots X_0
X_1 ]_{n_1} E_0 [X_1 X_0 \cdots ]_{n_2} \}$. Denote the space
spanned by $\{ [\cdots X_0 X_1 ]_{i} E_0 [X_1 X_0 \cdots ]_j | 0\leq
i,j\leq m-1 \}$  as $I^{'}_m$. We prove $I^{'}_m$  is stable under
the multiplication by $X_0 ,X_1$ from left side.

Case 1. $i<m-1$ and $i$ is odd, then the left most letter in
$[\cdots X_0 X_1]_i $ is $X_1$, so we write the word in discussion
as $ [X_1 \cdots X_0 X_1 ]_{i} E_0 [X_1 X_0 \cdots ]_j$. Since
$i<m-1$, $X_0 \cdot [X_1 \cdots X_0 X_1 ]_{i} E_0 [X_1 X_0 \cdots
]_j \in I^{'}_m $. And we have $X_1 \cdot [X_1 \cdots X_0 X_1 ]_{i}
E_0 [X_1 X_0 \cdots ]_j =$

$ X_1 ^2 [\cdots X_0 X_1 ]_{i-1} E_0 [X_1 X_0 \cdots ]_j = E_1
[\cdots X_0 X_1 ]_{i-1} E_0 [X_1 X_0 \cdots ]_j mod I^{'} _m$ by
$(3)$ of Definition 3.4. Using $(2)$ of Lemma 4.1, we have $E_1
[\cdots X_0 X_1 ]_{i-1} E_0 [X_1 X_0 \cdots ]_j =$

$ [\cdots X_0 X_1 ]_{m-1} E_0 [\cdots X_0 X_1 ]_{m-1}^{-1} [\cdots
X_0 X_1 ]_{i-1} E_0 [X_1 X_0 \cdots ]_j = \gamma [\cdots X_0 X_1
]_{m-1} E_0 [X_1 X_0 \cdots ]_j $ for some $\gamma \in \mathbb{C}$
by $(5)$ of Definition 3.4.

Case 2. $i<m-1$ and $i$ is even, then the word can be written as $
[X_0 \cdots X_0 X_1 ]_{i} E_0 [X_1 X_0 \cdots ]_j$. Since $i<m-1$,
$X_1 \cdot [X_0 \cdots X_0 X_1 ]_{i} E_0 [X_1 X_0 \cdots ]_j \in
I^{'}_m $. It is easy to see we have also $X_0 \cdot [X_0 \cdots X_0
X_1 ]_{i} E_0 [X_1 X_0 \cdots ]_j \in I^{'}_m $.

Case 3. $i=m-1$. It is easy to see $X_0  [\cdots X_0 X_1 ]_{m-1} E_0
[X_1 X_0 \cdots ]_j \in I^{'}_m $. And we have $X_1  [\cdots X_0 X_1
]_{m-1} E_0 [X_1 X_0 \cdots ]_j = X_1 E_1 [\cdots X_0 X_1 ]_{m-1}
[X_1 X_0 \cdots ]_j =$

$ l^{-1} E_1 [\cdots X_0 X_1 ]_{m-1} [X_1 X_0 \cdots ]_j =
l^{-1}[\cdots X_0 X_1 ]_{m-1} E_0 [X_1 X_0 \cdots ]_j \in I^{'}_m $.

Similarly we can prove $I^{'} _m$ is stable by multiplication by
$X_0 , X_1 $ from the right side. So $(1)$ is proved.

For $(2)$, denote the subspace in $B_{D_m}(\Pi ,\kappa)$  spanned by
$$ \{ [\cdots X_0 X_1 ]_i E_0 [X_1 X_0 \cdots ]_j | 0\leq i,j\leq
k-1 \} \cup \{ [\cdots X_1 X_0 ]_i E_1 [X_0 X_1 \cdots ]_j | 0\leq
i,j\leq k-1 \} $$ as $I^{'}_m$.   We only need to show $I_m =
I^{'}_m $ as well. Similar discussions by using Definition 3.3
(Notice the relation $(6)$ ) and Lemma 4.1 show that $I_m$ is
spanned by the set $\{ [\cdots X_0 X_1 ]_{n_1} E_0 [X_1 X_0 \cdots
]_{n_2} \} \cup \{ [\cdots X_1 X_0 ]_{n_1} E_1 [X_0 X_1 \cdots
]_{n_2} \}$. So it is enough to show that $I^{'}_m $ is stable under
the left and right multiplication by $X_0 ,X_1$. It suffice to show
the left case and the cases for words $ \{ [\cdots X_0 X_1 ]_i E_0
[X_1 X_0 \cdots ]_j | 0\leq i,j\leq k-1 \}$. Other cases are
similar.

Case 1. $i<k-1$ and $i$ is odd. Then the word in discussion can be
presented as $[X_1 \cdots X_0 X_1 ]_i E_0 [X_1 X_0 \cdots ]_j$.
Since $i<k-1$ we have $X_0 \cdot [X_1 \cdots X_0 X_1 ]_i E_0 [X_1
X_0 \cdots ]_j \in I^{'} _m $. On the other hand, $X_1 \cdot [X_1
\cdots X_0 X_1 ]_i E_0 [X_1 X_0 \cdots ]_j = X_1 ^2 [\cdots X_0 X_1
]_{i-1} E_0 [X_1 X_0 \cdots ]_j \in I^{'} _m$ by $(3),(6)$ of
Definition 3.3.

Case 2. $i<k-1$ and $i$ is even. The word can be written as $[X_0
\cdots X_0 X_1 ]_i E_0 [X_1 X_0 \cdots ]_j$. Since $i<k-1$ we have
$X_1 \cdot [X_0 \cdots X_0 X_1 ]_i E_0 [X_1 X_0 \cdots ]_j \in I^{'}
_m $. On the other hand,

$X_0 \cdot [X_0 \cdots X_0 X_1 ]_i E_0 [X_1 X_0 \cdots ]_j = X_0 ^2
[\cdots X_0 X_1 ]_{i-1} E_0 [X_1 X_0 \cdots ]_j \in I^{'} _m $ by
$(3),(5)$ of Definition 3.3.

Case 3. $i=k-1$.  If $k$ is odd, then the word can be written as
$[X_0 \cdots X_0 X_1 ]_{k-1} E_0 [X_1 X_0 \cdots ]_j$. Now we have

$X_0 \cdot [X_0 \cdots X_0 X_1 ]_{k-1} E_0 [X_1 X_0 \cdots ]_j \in
I^{'} _m$ by $(3),(5)$ of Definition 3.3. And we have

$X_1 \cdot [X_0 \cdots X_0 X_1 ]_{k-1} E_0 [X_1 X_0 \cdots ]_j =
\tau_0 ^{-1} \Phi ^0 _{m, \Pi}([\cdots X_0 X_1 ]_{m-1}) [\cdots X_1
X_0 ]_{k-1} ^{-1} E_0 [X_1 X_0 \cdots ]_j =$

$\tau_0 ^{-1} \Phi ^0 _{m, \Pi}([\cdots X_0 X_1 ]_{m-1}) [\cdots X_0
^{-1} X_1 ^{-1} ]_{k-1} E_0 [X_1 X_0 \cdots ]_j =  \gamma [\cdots
X_0 X_1 ]_{k-1}  E_0 [X_1 X_0 \cdots ]_j \mod I^{'} _m $ for some
$\gamma \in \mathbb{C}$. Where the first $"="$ in above equation is
by $(7)$ of Lemma 4.1, the third $"="$ is by iterated use of
Definition 3.3. The cases of $k$ being even and $i=k-1$ are easy to
see.

\end{proof}

The dihedral group $D_3$ is just the symmetric group $S_3$. We have
\begin{prop} The algebra $B_{D_3} (\Pi ,\kappa )$ in Definition 3.4 is
isomorphic the BMW algebra $B_3 (q,l)$.
\end{prop}
\begin{proof} By comparing Definition 1.1 for $B_3 (q,l) (B_{A_2}
(q,l))$ with Definition 2.1 for $B_{D_3} (\Pi ,\kappa )$, we only
need to prove  that the relation $5)$ of Definition 3.4 can be
induced from relations in Definition 1.1 for  $B_3 (q,l)$. First,
for  $B_3 (q,l)$, let $X$ be any word composed from $X_i ,E_i$
$(i=0,1)$, we claim there exist $\phi(X), \psi(X) \in
\mathbb{C}[q^{\pm}, l^{\pm}]$, such that $E_0 X E_0 =\phi(X)E_0$,
$E_1 X E_1 = \psi(X)E_1$. Denote the length of a word $X$ as $l(X)$.
Call letters $X_0 , E_0 $ $(X_1 , E_1)$ as of type 0 $(type 1 )$. By
relations in Definition 1.1, we only need to consider those $X$ in
which type 0 letters and type 1 letters appear alternatively and
$l(x)$ being odd. It is proved in Proposition 2.3 of \cite{CGW1}
that $E_0 E_1 E_0 = E_0$ and $E_1 E_0 E_1 = E_1$. So the claim is
true for $l(X)=0,1$. Suppose we have proved the claim for $l(X) \leq
2k+1$. Now suppose $l(X)= 2k+3$. For the case $E_0 X E_0$, we can
assume the letter $E_0$ doesn't appear in $X$ and $E_1$ appear in
$X$ at most once. Because otherwise we can use induction. So we only
need to consider those $X$ with $l(X)\leq 3$. Because if $l(X) \geq
5$, by assumption above $X= X_1 X_0 A_1 Y$ or $X= Y A_1 X_0 X_1$,
where $A_1 = X_1 \ or\ E_1$.  Now $E_0 X_1 X_0 A_1 Y E_0 = X_1 X_0
E_1 A_1 Y E_0 = \lambda X_1 X_0 E_1 Y E_0 = \lambda E_0 X_1 X_0 Y
E_0  $ for some $\lambda \in \mathbb{C}[q^{\pm}, l^{\pm}] $, which
implies $E_0 X_1 X_0 A_1 Y E_0 \in \mathbb{C}[q^{\pm}, l^{\pm}] E_0$
because $l ( X_1 X_0 Y )< l(X)$. So to prove the claim for $E_0 X
E_0$, we are left with the cases $X\in \{  X_1 X_0 X_1, E_1 X_0 X_1
, X_1 X_0 E_1 \} $, which can be certified by simple computations.
The cases for $E_1 X E_1$ can be proved similarly.

Since the Lawrence-Krammer representation factor through $B_3
(q,l)$, we see $\phi (X) =\Phi ^0 _{3,\Pi ,\kappa} (X) $ and $\psi
(X)= \Phi ^1 _{3,\Pi ,\kappa} (X)$, where $\Pi = \{ \alpha ,k  \}$
such that $q=\exp(\kappa \sqrt{-1} \pi k )$, $l= q^{\alpha -1}$.

\end{proof}

When $m$ is odd, suppose $(V,\rho)$ is a finite dimensional
representation of $Br_{D_m}( \Pi )$, where $\Pi =\{ k, \alpha  \}$.
Then on the trivial bundle $M_{D_m} \times V$ we have a
$D_m$-invariant flat connection
\begin{equation}
\rho(\Omega_{D_m}) = \kappa \sum _{i=0} ^{m-1} (k \rho(s_i) -
\rho(e_i)) \omega_i .
\end{equation}
It induces a flat connection $\rho(\bar{\Omega}_{D_m})$ on $M_{D_m}
\times_{D_m} V$, whose monodromy representation will be denoted as
$T_{\rho} : A_{D_m} \rightarrow GL(V)$. When $m=2k$ is even, suppose
$(V,\rho)$ is a finite dimensional representation of $Br_{D_m}( \Pi
)$, where $\Pi=\{ k_0 ,\alpha_0 ; k_1 ,\alpha_1 \}$. Similarly by
using the $D_m$-invariant flat connection
\begin{equation}
\rho(\Omega_{D_m}) = \kappa[  \sum _{i=0} ^{k-1} (k_0 \rho(s_{2i}) -
\rho(e_{2i})) \omega_{2i} +   \sum _{i=0} ^{k-1} (k_1 \rho(s_{2i+1})
- \rho(e_{2i+1})) \omega_{2i+1}] ,
\end{equation}
we obtain a monodromy representation $T_{\rho}$: $A_{D_m}
\rightarrow GL(V)$.

\begin{prop}(1) Suppose $m$ is odd.  Let the data $\Pi$ be generic so that $Br_{D_m}(
\Pi )$ is a semisimple algebra. Then for generic $\kappa$, the
representation $T_{\rho}$ of $A_{D_m}$ factors through the algebra
$B_{D_m}( \Pi ,\kappa )$. (2) Suppose $m$ is even. Let the data
$\Pi$ be generic so that $Br_{D_m}( \Pi )$ is a semisimple algebra.
Then for generic $\kappa$, the representation $T_{\rho}$ of
$A_{D_m}$ factors through the algebra $B_{D_m}(\Pi ,\kappa )$ .

\end{prop}

\begin{proof}
We prove $(2)$, the proof of $(1)$ is similar and easier. So let $m$
be even. Denote the natural quotient map from $Br_{D_m}(\Pi)$ to
$\mathbb{C}D_m$ as $\pi_m$. Suppose $\rho_1 , \cdots ,\rho_v$ are
all the irreducible representations of $\mathbb{C}D_m$, denote the
representation of $Br_{D_m}(\Pi)$ induced from $\rho_i$ through
$\pi_m$ as $\bar{\rho}_i$. As in Lemma 2.4, we have two more
representation $\rho_{LK,c_0}, \rho_{LK, c_1}$ of $B_{D_m}(\Pi)$. By
Theorem 6.1 ,Theorem 6.3 of \cite{CH}, we know for generic data
$\Pi$, the algebra $Br_{D_m}(\Pi)$ is semisimple and $\bar{\rho}_1
,\cdots ,\bar{\rho}_v ,\rho_{LK,c_0}, \rho_{LK, c_1} $ are all the
irreducible representations of $Br_{D_m}(\Pi)$. So because of
semisimplicity of $Br_{D_m}(\Pi)$, we only need to prove the
proposition when $\rho$ is one of above listed irreducible
representations. Set $q_i = \exp (\sqrt{-1} \pi \kappa k_i )$, $l_i
= \exp (\kappa \alpha_i \pi \sqrt{-1}) / q_i$. Set $e^{\rho} _i =
\frac{l_i}{ q^{-1}_i -q_i } (T_{\rho} (\sigma_i ) - q_i ) (T_{\rho}
(\sigma_i) + q^{-1}_i)$ for $i=0,1$. We want to show the map
$\phi_m:  E_i \mapsto e^{\rho} _i $, $X_i \mapsto T_{\rho}
(\sigma_i)$ extends to an algebraic morphism. Recall by Lemma 2.2,
$T_{\rho} (\sigma_i)$ is conjugated to $\rho (s_i) \exp (\kappa
\sqrt{-1} \pi (k_i \rho (s_i) - \rho (e_i) ))$ where $k_i = k_0 $
for even $i$, and $k_i = k_1$ for odd $i$. Suppose $\rho =
\bar{\rho}_l$, then $k_i \rho(s_i) - \rho (e_i ) = k_i \rho(s_i)$
can be presented by a diagonal matrix whose diagonal elements being
in $\{ \pm k_i \}$. So we know $e^{\rho} _i =0$ and $\phi_m$ keep
all relations in Definition 3.3 and these cases are done. Now
suppose $\rho = \rho_{LK,c_0}$.  By $(1)$ of lemma 3.2, we see
$\phi_m$ satisfy $(4)$ of Definition 3.3. The cases for
$(1),(2),(3)$ of Definition 3.3 are trivial. Since $T_{\rho}$ is a
subrepresentation of $T_{LK}$, by definition of $\Phi^i _{m, \Pi}
(X)$ we know relation $(5)$ of Definition 3.3 is also satisfied by
$\phi_m$. Since $T_{\rho_{LK,c_0}} (e_1 )=0$, we see $e^{\rho} _1
=0$ so relation $(6)$ of Definition 3.3 is also satisfied by
$\phi_m$.  The case for $\rho = \rho_{LK, c_1}$ is similar.
\end{proof}

\begin{rem}
If $(V,\rho )$ is the infinitesimal Lawrence-Krammer representation
of $Br_m (\Pi)$ $($in Lemma 2.3$)$ , then the connection
$\rho(\Omega_{D_m})$ of equation $(7)$ $($ $(8)$ $)$ coincide with
the connection $\Omega_m$ of equation $(1)$ $($ $(2)$ $ )$. So above
proposition shows the generalized Lawrence-Krammer representation of
$A_{D_m}$ factor through the algebra $B_{D_m} (\Pi ,\kappa )$. It
thus provide us a way to compute the monodromy of Marin's flat
connections by using the generalized BMW type algebras, at least in
cases of dihedral groups.
\end{rem}

\begin{thm} (1) When $m$ is odd, for generic data $\Pi ,\kappa$, the algebra
$B_{D_m}(\Pi ,\kappa )$ is semisimple with dimension $2m+ m^2$; (2)
When $m$ is even, for generic data $\Pi$,  the algebra $B_{D_m} (\Pi
,\kappa )$ is semisimple with dimension $2m + \frac{m^2}{2}$.
\end{thm}
\begin{proof}
First we prove $(1)$. According to Proposition 4.1,  denote the
quotient map from $B_{D_m}(\Pi ,\kappa )$ to the Hecke algebra
$H_{D_m}(q)$ as $\pi$. For generic $q$, $H_{D_m}(q)$ is semisimple
and we denote its irreducible representations as $\rho_1 ,\cdots
,\rho_k $. Through $\pi$, each $\rho_i$ induces an irreducible
representation $\bar{\rho}_i$ of $B_{D_m}(\Pi ,\kappa )$. These
representations are different  from each other and we have
$\sum_{i=1} ^k (\dim \bar{\rho} _i )^2 = \dim H_{D_m}(q)=2m. $
Because for any $1\leq i\leq k$, the annihilating polynomial of
$\bar{\rho}_i (X_j ) $ has degree $2$, and the annihilating
polynomial of $\bar{T}_{LK}(X_j )$ has degree $3$, so by Lemma 3.5
we know $\bar{T}_{LK}$ is an irreducible representation different
with any $\bar{\rho}_i$. So by Wedderburn-Artin theorem, we have
$\dim B_{D_m} (\Pi ,\kappa ) \geq \sum_{i=1} ^k (\dim \bar{\rho} _i
)^2 + (\dim \bar{T}_{LK})^2 = 2m+ m^2$. So by $(1)$ of Lemma 4.2 and
by Wedderburn-Artin theorem again, we know for generic data $\Pi
,\kappa$, $B_{D_m} (\Pi ,\kappa )$ is a semisimple algebra with
dimension $ 2m+ m^2$. The proof of $(2)$ is similar, by using $(2)$
of Lemma 4.2 and the Wedderburn-Artin theorem.

\end{proof}

\begin{rem} By Theorem 6.1 of $\cite{CH}$, we have $\dim Br_{D_m}(\Upsilon)= m^2
+2m$ for odd $m$, and $\dim Br_{D_m}(\Upsilon)= m^2 +\frac{m}{2}$
for even $m$. Above Theorem 4.1 and Proposition 4.2, 4.3 say that
the algebra $B_{D_m}(\Pi ,\kappa )$ is a satisfactory deformation of
the Brauer type algebra $Br_{D_m}(\Upsilon).$
\end{rem}
Now let $W_{\Gamma}$ be a  Coxeter group with Coxeter matrix $\Gamma
= (m_{i,j})_{n\times n}$. Which has a presentation as in Theorem
2.2. For each $1\leq i\leq n$, choose a pair of constants $(k_i  ,
\alpha_i ) \in \mathbb{C}^{\times} \times \mathbb{C}$ such that
$(k_i  , \alpha_i )=(k_j  , \alpha_j ) $ if $s_i$ is conjugated to
$s_j$. Denote the data consisting of $\{ (k_i ,\alpha_i) \}_{i=1}
^m$ as one symbol $\Pi$. Let $\kappa \in \mathbb{C}$. Set $q_i =
\exp(\kappa \sqrt{-1}\pi k_i )$, $l_i = (q_i )^{\alpha_i -1}$ for
$1\leq i\leq n$. The definition of BMW type algebras of Dihedral
groups naturally inspire us to present the following definition of
$B_{W}(\Pi, \kappa)$, the BMW type algebra of type $W$. These
algebras will be generated by elements $X_i ,E_i$, $1\leq i\leq n$.
For $1\leq i <j\leq n$, denote the parabolic subgroup of $W$
generated by $s_i ,s_j$ as $W_{i,j}$. Then $W_{i,j}$ is isomorphic
to $D_{m_{i,j}}$. Denote the data $\{ (k_i ,\alpha_i), ( k_j ,
\alpha_j ) \}$ as $\Pi_{i,j}$. Denote the free monoid generated by
$X_i , E_i , X_j , E_j$ as $\Lambda _{i,j}$. Repeat the definition
of function $\Phi^i _{m, \Pi ,\kappa }$  in Definition 3.1, but
replace the roles of $( D_m , X_0 , E_0 , X_1 , E_1 , \Lambda , \Pi
)$ by $( W_{i,j} , X_i , E_i , X_j , E_j , \Lambda_{i,j} ,\Pi_{i,j}
)$, we obtain functions $\Phi^i _{m_{i,j} , \Pi_{i,j}} $ and $\Phi^j
_{m_{i,j} , \Pi_{i,j}}$ from $\Lambda_{i,j}$ to $\mathbb{C}$, which
are similar to $\Phi^0 _{m, \Pi ,\kappa}$ and $\Phi^1 _{m, \Pi
,\kappa}$ in Definition 3.1 respectively.
\begin{defi} The algebra $B_{\Gamma}(\Pi ,\kappa)$ is generated by elements $X^{\pm} _i
,E_i$, $1\leq i\leq n$ with the following relations.

$(1)$. $X_i X^{-}_i = X^{-}_i X_i =1$ for any $i$;

$(2)$. $[X_i X_j \cdots ]_{m_{i,j}} = [X_j X_i \cdots ]_{m_{i,j}}$
for any $i\neq j$;

$(3)$. $(q_i ^{-1} -q_i ) E_i = l_i (X_i -q_i )(X_i + q_i ^{-1} )$
for any $i$;

$(4)$. $X_i E_i = E_i X_i = l_i ^{-1} E_i $ for any $i$;

$(5)$. $E_i X E_i = \Phi^i _{m_{i,j}, \Pi_{i,j},\kappa}(X)E_i $ and
$E_j X E_j = \Phi^j _{m_{i,j}, \Pi_{i,j},\kappa}(X)E_j $ for any
$i<j$ and $X\in \Lambda_{i,j}$;

$(6)$. $E_i X E_j =0$ for any $X\in \Lambda_{i,j}$ if $m_{i,j} $ is
an even number greater than $2$.

\end{defi}

\section{ General BMW Type Algebras}
Let $\Gamma =(m_{i,j})_{n\times n}$ be a Coxeter matrix. For each
$1\leq i\leq n$, choose $q_i \in \mathbb{C}^{\times}$ such that $q_i
=q_j$ if $s_i $ is conjugated to $s_j$. Denote the data $\{ q_1 ,
\cdots , q_n \}$ as $\bar{q}$. The type $\Gamma$ Hecke algebra
$H_{\Gamma} (\bar{q})$ is the quotient algebra of the group algebra
$\mathbb{C} A_{\Gamma}$ to the ideal generated by $\{ (\sigma_i -q_i
)(\sigma_i + q^{-1}_i ) \}$. The following Proposition is evident.

\begin{prop} The map $X_i \mapsto \sigma_i $, $E_i \mapsto 0$ $(1\leq i\leq
n)$ extends to a surjection from $B_{\Gamma}(\Pi ,\kappa)$ to
$H_{\Gamma} (\bar{q})$ if the data $\Pi$ and $\bar{q}$ satisfies
$q_i = \exp \kappa \sqrt{-1} \pi k_i $ for $1\leq i\leq n$.
\end{prop}

The following Theorem 5.1 is proved using similar method as in the
proof of Proposition 2.9 of \cite{CGW1}. Suppose $\Gamma
=(m_{i,j})_{n\times n}$ is a finite type Coxeter matrix. Let $w_0$
be the longest element of $W_{\Gamma}$ and denote the length of
$w_0$ as $L$. A sequence  in $[1,n]$ is a sequence $(i_1 , i_2
,\cdots , i_k )$ where $1\leq i_v \leq n$ for any $v$. The sequence
$(i_1 , i_2 ,\cdots , i_k )$ is called reducible if $i_v = i_{v+1}$
for some $v$. Otherwise it is called irreducible. A basic
transformation of a sequence $Q$ is, replacing  a subsequence  $[i,j
, \cdots ]_{m_{i,j}}$ in $Q$ by $[j,i,\cdots ]_{m_{i,j}}$. For
example, transform $(3,1,2,1,2,4)$ to $(3,2,1,2,1,4)$ when
$m_{1,2}=4$. By the theory of Coxeter groups, we see the following
fact:  if a sequence $Q$ has length greater than $L$, then $Q$ can
be transformed into a reducible sequence.

\begin{thm} If $\Gamma$ is of finite type, then $B_{\Gamma}(\Pi ,\kappa)$ is a finite dimensional algebra.
\end{thm}
\begin{proof}   Denote the set of words composed by $\{ X_i  ,E_i \}_{1\leq
i\leq n}$ as $S$. By $(3)$ of Definition 4.1, we see $B_{\Gamma}(\Pi
,\kappa )$ can be spanned by $S$. Let $w_0$ be the longest element
of $W_{\Gamma}$, denote the length of $w_0$ as $L$.  For a word $W$,
denote the length of $W$ as $l(W)$. Let $S_0 =\{ W\in S | l(W)\leq L
\}$. Denote the subspace of $B_{\Gamma}(\Pi ,\kappa)$ spanned by
$S_0 $ as $A$.  The theorem follows from the following assertion:
Any word in $S$ lie in $A$ as an element of the algebra. Let $W\in
S$. When $l(W)<L$, the assertion is evident. Suppose we have prove
that for any $W\in S$ with $l(W)\leq K$,  we have $W\in A$. Here we
can ask $K\geq L$. Let $W\in S$ is a word with length $K+1$. Each
word $W\in S$ determine a sequence of indices $Q_W$. For example, if
$W= E_1 X_2 X_3 E_2 $ then  $Q_W =(1,2,3,2)$. Since the length of
$Q_W$ is greater than $L$, we see $Q_W$ can be transformed into a
reducible sequence by $M$ times of basic transformations.  Then we
do induction on $M$. If $M=0$, then $Q_W$ is reducible, which mean
$W$ contains a subword with one of the following type: $E_i E_i$,
$X_i E_i $, $E_i X_i$, $X_i X_i$.  By relation $(3)$, $(4)$ of
Definition 4.1 and by induction we see $W\in A$.  Now suppose  the
assertion of $l(W)\leq K$ and $l(W)=K+1 , M\leq J$ are proved.
Suppose we have a word $W$ such that $Q_W$ is irreducible,
$l(W)=K+1$ and $M= J+1$. Then there is a series of sequence $Q_1 ,
Q_2 ,\cdots ,Q_{J+2}$ such that $Q_1 = Q_W$, $Q_{J+2}$ is reducible
and $Q_{i+1} $ is obtained by doing one basic transformation to
$Q_i$. Denote the subword of $W$ in the position where the baisc
transformation to $Q_2$ occur as $V$. And suppose the basic
transformation from $Q_1$ to $Q_2$ is by replacing a subsequence
$[i,j\cdots ]_{m_{i,j}}$ with $[j,i\cdots ]_{m_{i,j}}$. For example
when $W= X_1 E_2 X_3 X_2 E_3$, and $Q_2 = (1,3,2,3,2)$ then $V= E_2
X_3 X_2 E_3$.  First we observe if there are two $E_i$, or two $E_j$
in $V$ then we can shorten $V $ by Relation $(5)$ of Definition 4.1
and the assertion for $W$ can be proved. So we can suppose there is
at most one $E_i $ and one $E_j$ in $V$.

Case 1. $m_{i,j}=2k+1$ is odd, $V= A E_i E_j C$ where $A, C$ contain
only $X_i ,X_j$. Denote the length of $A,C$ as $a,c$ respectively,
so by above arguments $A=[\cdots X_i X_j]_a$ and $C=[\cdots X_j X_i
]_c $.  We have $a+c=2k-1$. Suppose $a>c$, so $a\geq k$. By applying
$(2)$ of Lemma 4.1 to the subalgebra generated by $\{ X_i , X_j ,
E_i ,E_j \}$, we have $A E_i E_j C = [\cdots ]_{2k-a} E_j D E_j C =
\lambda [\cdots]_{2k-a} E_j C$, where $[\cdots]_{2k-a}$ denote some
word with length $2k-a$, and $D$ is some word. Since
$l([\cdots]_{2k-a} E_j C)\leq K$, by replace $X^{-1} _i$ in
$[\cdots]_{2k-a}$ with linear sum of $X_i , E_i $, we complete the
proof of this case by induction. The case of $a<c$ is similar.

Case 2. $m_{i,j}=2k+1$ is odd, $V= A E_i B E_j C $ Where $A,B,C$ are
nonempty words contain only $X_i ,X_j $. Denote the length of
$A,B,C$ as $a,b,c$ respectively. We have $a\neq c$ because otherwise
 $V= A E_i B E_i C$ which contradicts our assumption. We can
suppose $a>c$. The cases of $a<c$ is similar.   Again by using $(2)$
of Lemma 4.1, $ A E_i B E_j C = [\cdots ]_{2k-a} E_j D E_j C
=\lambda [\cdots ]_{2k-a} E_j C$. Now $l (\lambda [\cdots ]_{2k-a}
E_j C)= 2k-a +1+c < 2k-1 $ so by induction we complete the proof of
these cases.

Case 3. $m_{i,j}=2k+1$ is odd, $V= A  E_i B $ where $A,B$ are
nonempty words contain only $X_i ,X_j $. Denote the length of $A,B$
as $a,b$ respectively. Suppose $B$ is nonempty. Relation $(3)$ of
Definition 4.1 gives us an identity $X_i = \lambda X^{-1}_i +\mu E_i
 +\gamma$. Replace every $X_i$ in  $B$ with $ \lambda X^{-1}_i +\mu E_i
+\gamma$ and every $X_j $ in B with $\lambda X^{-1}_j +\mu E_j
+\gamma$, we see $V= \lambda^{'} A E_i [X^{-1}_j X^{-1}_i \cdots ]_b
+ U $ where $U$ is a liner sum of words containing two $E_i$, or
containing more than one $E_i$ and containing one $E_j$ with length
$K+1$, or with smaller length.  By Case 1 and 2, we see $U\in A$.
And by $(2)$ of Lemma 4.1, we have $A E_i [X^{-1}_j X^{-1}_i \cdots
]_b = [\cdots X^{-1}_j X^{-1}_i]_b E_j A^{'}$. By replacing
$X^{-1}_i , X^{-1}_j$ in  the subword $[\cdots X^{-1}_j X^{-1}_i]_b$
we see $W= $ a linear sum of words with sequence $Q_2 $ $\mod$ $A$.
Since $Q_2 $ can be transformed to an reducible sequence with fewer
basic transformations, we proved these cases.

Case 4. $m_{i,j}=2k+1$ is odd, $V= [X_i X_j \cdots ]_{2k+1}$.  By
replacing the subword $V$ in $W$ by $[X_j X_i \cdots ]_{2k+1}$, we
see the resulted word (equal $W$ ) has sequence $Q_2$, then we use
induction.

For the cases when $m_{i,j}= 2k$, the same argument shows we only
need to consider those cases when $V$ contain at most one $E_i$ and
one $E_j$. Now $(6)$ of Definition 4.1 show we only need to consider
the cases  (1) that $V$ contain one $E_i$ or one $E_j$; (2) the case
$V$ contain only $ X_i ,X_j  $. We can prove  the cases $(1)$
similarly as above case 3, and prove the cases $(2)$ as above case
4.
\end{proof}
When $\Gamma$ is a simply laced type Coxeter group, then all
reflections of $W_{\Gamma}$  lie in the same conjugacy class. So the
data $\Pi$ of $B_{\Gamma}(\Pi ,\kappa)$ consists of $\{ k, \alpha
\}$ essentially.  As a corollary of Proposition 4.2, we have
\begin{prop} The algebra $B_{\Gamma} (\Pi ,\kappa)$ defined in Definition 4.1 is isomorphic to the simply laced BMW
algebra  $B_{\Gamma} (q,l)$ of $\cite{CGW1}$ for $q=\exp(\kappa
\sqrt{-1}\pi k) $, $l=\frac{exp(\alpha \pi \sqrt{-1})}{q}$.
\end{prop}

 Let $W$ be a finite type Coxeter group. Let $R$ be the set of
reflections in $W$. Suppose $(V, \rho)$ is a finite dimensional
representation of $Br_W (\Pi )$. Where $\Pi = \{
 k_s ,\alpha_s \}_{s\in R}$.  Then by using the KZ connection ( Definition 2.3), we have a flat, $W$-invariant connection
 $\rho(\Omega_W)= \kappa \sum_{s\in R} (k_s \rho (s) - \rho(e_s) ) \omega_s $ on $M_W \times V$, which induces a representation $T_{\rho} : A_{W}\rightarrow
 GL(V)$. Suppose $W$ is of rank $n$ and has a canonical set of generators $\{ s_i \}_{1\leq i\leq n}$ as in Theorem 2.2. Suppose the Artin group $A_W$ has a canonical set
 of generators $\{ \sigma_i \}_{1\leq i\leq n}$ as in Theorem 2.3.
 \begin{thm}
For generic $\Pi$ and $\kappa$, the representation $T_{\rho}$ factor
through $B_W ( \Pi ,\kappa )$.
 \end{thm}
 \begin{proof}
We want to show if set $X_i = T_{\rho}(\sigma_i )$ and $E_i  =
\frac{l_i }{q_i ^{-1}-q_i } (T_{\rho}(\sigma_i ) - q_i )
(T_{\rho}(\sigma_i) + q_i ^{-1}) $ for $1\leq i \leq n$, where $q_i
= \exp(\kappa k_{s_i} \pi \sqrt{-1})$ and $l_i = \exp(\kappa
\alpha_i \pi \sqrt{-1})/q_i$, then $\{ X_i , E_i \}_{1\leq i\leq n}$
satisfies the relations in Definition 4.1. If $W$ is a dihedral
group $D_m$, this is true by Proposition 4.3. The cases for general
$W$ can be reduced to the cases for dihedral groups by using Theorem
2.4 as follows. Recall the notations " $\Delta $, $H_{s_i}$" above
Theorem 2.2, and the notations "$\Pi$, $\Pi_{i,j}$ " before
Definition 4.1. For $1\leq i< j\leq n$, first we suppose $m_{i,j}$
is odd.  Denote $L= H_{s_i} \cap H_{s_j}$, define $R_L $, $W_{L}$,
$A_{W_L}$ and $M_{W_L}$ as in section 2 (above Theorem 2.4). We can
identify $A_{W_L}$ with the Artin group $A_{D_{m_{i,j}}}$ in
suitable way, so the morphism $\bar{\lambda}_L$ from $A_{W_L}$ to
$A_W$ maps $\sigma_0 , \sigma_1$ to $\sigma_i , \sigma_j$
respectively.

In this case we have
 a natural isomorphism $\phi_{i,j}$ from $W_{D_{m_{i,j}}}$ to $W_{L}$ extending
 the map $s_0 \mapsto s_i $, $s_1 \mapsto s_j$, and we can identify
 $M_{W_{D_{m_{i,j}}}}$ with $M_{W_L}$ ( actually $ M_{W_L}\cong M_{W_{D_{m_{i,j}}}}\times \mathbb{C}^{n-2}$, but it will not
 make any difference. ) The morphism $\phi_{i,j}$ can be extended to
 a morphism from $Br_{D_{m_{i,j}}} (\Pi_{i,j}) $ to $Br_W
 (\Pi)$ by sending $e_0 , e_1 \in Br_{D_{m_{i,j}}} (\Pi_{i,j})$ to
 $e_i ,e_j \in Br_W
 (\Pi)$ respectively. And we know the morphism
 $\bar{\lambda}_{L}$ maps $A_{W_L}$ isomorphically to the parabolic
 subgroup of $A_W$ generated by $\sigma_i , \sigma_j$.  Through $\phi_{i,j}$ we have a representation $\rho \circ  \phi_{i,j} :Br_{D_{m_{i,j}}} (\Pi_{i,j})\rightarrow gl(V)
 $. As in section 2, set $\Omega^{'}=\kappa \sum_{s\in R_L} ( k_s  \rho (s)- \rho (e_s )  ) \omega_s
 $. It isn't hard to see $\Omega^{'}$ coincide with the connection
 $ \Omega_{\rho \circ  \phi_{i,j} }$ on $M_{W_L}$. So by Proposition
 4.3, if we set  $X_0 = T_L (\sigma_0)$ and $X_1 = T_{L}(\sigma_1)$,
 then $X_0 , X_1 $ satisfy the relations in Definition 3.4 for
 $m=m_{i,j}$ and suitable data $\Pi$. ( since by $(3)$ of Definition 3.4, $E_i$ can be presented by $X_i$ for generic data, all those relations
 can be seen as relations between $X_0 ,X_1$. ) By Theorem 2.4,
 there is a element $A\in GL(V)$, so $T\circ \bar{\lambda}_L (\rho_i) = A T_L (\rho_i )
 A^{-1}$ for $i=0,1$. Since $T_{\rho}\circ \bar{\lambda}_L (\rho_0)= T_{\rho}(\sigma_i ), T_{\rho}\circ \bar{\lambda}_L
 (\rho_1)=T(\sigma_j)$, we have that $X_l = T_{\rho}(\sigma_l )$ and $E_l  =
\frac{l_l }{q_l ^{-1}-q_l } (T_{\rho}(\sigma_l ) - q_l )
(T_{\rho}(\sigma_l) + q_l ^{-1}) $ $(l=i,j)$ satisfy the relation
$(5)$ of Definition 4.1.

Similarly we can prove if $m_{i,j}$ is even, then $X_i , X_j ,E_i
,E_j$ satisfy relation $(5), (6)$ of Definition 4.1, and we can
prove $X_i , E_i $ satisfy relation $(3),(4)$, by using Theorem 2.4.

 \end{proof}


\end{document}